\documentclass[11pt]{amsart}
\usepackage{amssymb,latexsym}

\newdimen\AAdi%
\newbox\AAbo%
\def\AAk#1#2{\s_etbox\AAbo=\hbox{#2}\AAdi=\wd\AAbo\kern#1\AAdi{}}%
\def\AAr#1#2#3{\s_etbox\AAbo=\hbox{#2}\AAdi=\ht\AAbo\raise#1\AAdi\hbox{#3}}%
\font\tenmsb=msbm10 at 12pt
\font\sevenmsb=msbm7 at 8pt
\font\fivemsb=msbm5 at 6pt
\newfam\msbfam
\textfont\msbfam=\tenmsb
\scriptfont\msbfam=\sevenmsb
\scriptscriptfont\msbfam=\fivemsb
\def\Bbb#1{{\tenmsb\fam\msbfam#1}}

\newcommand{\beq}{\begin{equation}}
\newcommand{\eeq}{\end{equation}}
\newcommand{\ba}{\begin{array}}
\newcommand{\ea}{\end{array}}

\begin{document}

\newtheorem{thm}{Theorem}
\newtheorem{lem}{Lemma}
\newtheorem{cor}{Corollary}
\newtheorem{rem}{Remark}
\newtheorem{pro}{Proposition}
\newtheorem{defi}{Definition}
\newcommand{\noi}{\noindent}
\newcommand{\dis}{\displaystyle}
\newcommand{\mint}{-\!\!\!\!\!\!\int}

\def \bx{\hspace{2.5mm}\rule{2.5mm}{2.5mm}} \def \vs{\vspace*{0.2cm}}
\def\hs{\hspace*{0.6cm}}
\def \ds{\displaystyle}
\def \p{\partial}
\def \O{\Omega}
\def \o{\omega}
\def \b{\beta}
\def \m{\mu}
\def \l{\lambda}
\def\L{\Lambda}
\def \ul{u_\lambda}
\def \D{\Delta}
\def \d{\delta}
\def \s{\sigma}
\def \e{\varepsilon}
\def \a{\alpha}
\def \tf{\tilde{f}}
\def\cqfd{%
\mbox{ }%
\nolinebreak%
\hfill%
\rule{2mm} {2mm}%
\medbreak%
\par%
}
\def \pr {\noindent {\it Proof.} }
\def \rmk {\noindent {\it Remark} }
\def \esp {\hspace{4mm}}
\def \dsp {\hspace{2mm}}
\def \ssp {\hspace{1mm}}

\def \u{u_+^{p^*}}
\def \ui{(u_+)^{p^*+1}}
\def \ul{(u^k)_+^{p^*}}
\def \energy{\int_{\R^n}\u }
\def \sk{\s_k}
\def \mo{\mu_k}
\def\cal{\mathcal}
\def \I{{\cal I}}
\def \J{{\cal J}}
\def \K{{\cal K}}
\def \OM{\overline{M}}

\def\fk{{{\cal F}}_k}
\def\M1{{{\cal M}}_1}
\def\Fk{{\cal F}_k}
\def\Fl{{\cal F}_l}
\def\FF{\cal F}
\def\Gk{{\Gamma_k^+}}
\def\n{\nabla}
\def\uuu{{\n ^2 u+du\otimes du-\frac {|\n u|^2} 2 g_0+S_{g_0}}}
\def\uuug{{\n ^2 u+du\otimes du-\frac {|\n u|^2} 2 g+S_{g}}}
\def\sku{\sk\left(\uuu\right)}
\def\qed{\cqfd}
\def\vvv{{\frac{\n ^2 v} v -\frac {|\n v|^2} {2v^2} g_0+S_{g_0}}}
\def\vvs{{\frac{\n ^2 \tilde v} {\tilde v}
 -\frac {|\n \tilde v|^2} {2\tilde v^2} g_{S^n}+S_{g_{S^n}}}}
\def\skv{\sk\left(\vvv\right)}
\def\tr{\hbox{tr}}
\def\pO{\partial \Omega}
\def\dist{\hbox{dist}}
\def\RR{\Bbb R}\def\R{\Bbb R}
\def\C{\Bbb C}
\def\B{\Bbb B}
\def\N{\Bbb N}
\def\Q{\Bbb Q}
\def\Z{\Bbb Z}
\def\PP{\Bbb P}
\def\EE{\Bbb E}
\def\F{\Bbb F}
\def\G{\Bbb G}
\def\H{\Bbb H}
\def\SS{\Bbb S}\def\S{\Bbb S}

\def\lcf{{locally conformally flat} }

\def\circledwedge{\setbox0=\hbox{$\bigcirc$}\relax \mathbin {\hbox
to0pt{\raise.5pt\hbox to\wd0{\hfil $\wedge$\hfil}\hss}\box0 }}

\date{}
\title[On a fully nonlinear Yamabe problem]
{On a fully nonlinear Yamabe problem}
\author{Yuxin Ge}
\address{Laboratoire  d'Analyse et de Math\'ematiques Appliqu\'ees,
CNRS UMR 8050,
D\'epartement de Math\'ematiques,
Universit\'e Paris XII-Val de Marne, \\61 avenue du G\'en\'eral de Gaulle,
94010 Cr\'eteil Cedex, France}
\email{ge@univ-paris12.fr}
\author{Guofang Wang}
\address{Max-Planck-Institute for Mathematics in
the Sciences\\ Inselstr. 22-26, 04103 Leipzig, Germany}
\email{gwang@mis.mpg.de}
\begin{abstract} We solve the $\s_2$-Yamabe problem for a non \lcf manifold of dimension
$n>8$.
\end{abstract}

\maketitle

{\bf R\'esum\'e :} On r\'esout le probl\`eme de $\sigma_2$-Yamabe pour des vari\'et\'es riemanniennes compactes sans
bord non localement conform\'ement plates de dimension $n>8$.\\

\quad \quad  Dedicated to Professor W. Y. Ding on the occasion of his 
60th birthday

\

\

\section{Introduction}
Let $(M,g_0)$ be a compact, oriented Riemannian manifold
with metric $g_0$ and $[g_0]$ the conformal class of $g_0$. Let 
$Ric_g$ and $R_g$ be the Ricci tensor and scalar curvature of
$g$ respectively.
The Schouten tensor of the metric $g$ is defined by
\[S_g=\frac 1{n-2}\left(Ric_g-\frac {R_g}{2(n-1)}\cdot g\right).\]
The Schouten tensor plays an important role in conformal geometry. 
 Let $\s_k$ be the $k$th elementary symmetric function. For a symmetric
 $n\times n$ matrix $A$, set $\s_k(A)=\s_k(\L)$, where $\L=(\l_1,\l_2,\cdots, \l_n)$
 is the set of eigenvalues of $A$.
The  {\it $\s_k$-scalar curvature} 
of $g$ is defined by
\[\s_k(g):=\s_k (g^{-1}\cdot S_g),\]
where $g^{-1}\cdot S_g$ is locally defined by $(g^{-1}\cdot S_g)^i_j=
\sum_k g^{ik}(S_g)_{kj}$, see \cite{Jeff1}. Note that $\s_1(g)=\frac 1{2(n-1)}R_g$. It is an interesting question to
find a metric  $g$ in a given conformal class $[g_0]$ such that
\beq\label{0.1}
\s_k(g)=c
\eeq
for some constant $c$. Since  the Schouten tensors $S_g$ and $S_{g_0}$
of conformal metrics $g=e^{-2u}g_0$
and $g_0$ have the following relation
\[S_g=\n ^2 u+du\otimes du-\frac {|\n u|^2} 2 g_0+S_{g_0},\]
Equation (\ref{0.1}) is equivalent to the following fully nonlinear equation
\beq\label{0.2}
\s_k\left(\uuu\right)=c e^{-2ku},\eeq
for some constant $c$. When $k=1$, it is the well-known Yamabe equation.

Let
\[\Gamma_k^+=\{\Lambda=(\l_1,\l_2,\cdots, \l_n)\in \R^n\,|\,
\sigma_j(\Lambda)>0, \forall j\le k\}\]
 be Garding's cone.  A metric $g$ is said to be $k$-positive
 or simply $g\in \Gamma_k^+$
if $\s_j(g)(x)>0$ for any $j\le k$ and at every point $x\in M$. 
If $g=e^{-2u}g_0$,
we say $u$ is $k$-admissible if $g$ is $k$-positive. 
In this paper we consider the following

\

\noindent{\it $\s_k$-Yamabe problem:} Let $g_0\in \Gamma_k^+$. Find a conformal metric
$g\in [g_0]\cap \Gamma_k^+$ such that 
\[\s_k(g)=c\]
for some constant $c$.\\

\

The study of  the fully nonlinear Equations (\ref{0.1}) was initiated by Viaclovsky.
Since then there is a lot of work concerning these equations. Here, we just mention
some results directly related to the existence of the $\s_k$-Yamabe problem.
This problem has been solved in the following cases. When $k=n$, under a sufficient
condition, Viaclovsky proved the existence in \cite{Jeff4}. 
When $n=2k=4$, which is an important case,
Chang-Gursky-Yang solved the problem in \cite{CGY2}. 
See also \cite{CGY1} and \cite{GV4}.
When the underlying manifold is locally conformally flat, this problem was solved
by Guan-Wang \cite{GW2} and Li-Li \cite{LiLi} independently. See also
\cite{BV}. Note that when the underlying manifold $(M,g_0)$
is locally conformally flat and $g\in \Gamma_k^+$ with $k\ge n/2$,  the universal cover of $M$ is
conformally equivalent to a spherical space form \cite{GVW}.
When $k>n/2$, the $\s_k$-Yamabe problem was solved
by Gursky-Viaclovsky in \cite{GV2}.
See also their earlier work \cite{GV1}.

In this paper,
we consider the case $k=2$. In this case, Equation (\ref{0.2}) is a variational problem,
which was observed by Viaclovsky in \cite{Jeff1}. This is crucial for our method 
presented here. Our main result in this paper is

\begin{thm}\label{maintheorem}
Let $(M^n,g_0)$ be a compact, oriented Riemannian manifold with $g_0\in \Gamma_2^+$.
When $n>8$ and  the Weyl tensor $W_{g_0}\not =0$,
then there is a conformal metric $g\in [g_0]\cap \Gamma_2^+$ such
that 
\[\s_2(g) =c\]
for some constant $c$.
\end{thm}

Combining the results of \cite{GW2} and \cite{LiLi}, the $\s_2$-Yamabe problem is solvable 
if $n>8$. Like the ordinary
Yamabe problem, there is a well-known difficulty --the loss of compactness
of Equation (\ref{0.1}). Another more difficult problem is the fully nonlinearity of (\ref{0.1}).
Our result here is an analogue of the result of Aubin \cite{Aubin} 
for the ordinary Yamabe problem.
 Even the ideas of proof
are quite similar. However the techniques to realize these  ideas
become more delicate due to the fully nonlinearity.

Set ${\cal C}_2=\{g\in [g_0] \,|\, g\in\Gamma_2^+\}$ and define a
Yamabe type constant by
\[Y_2(M, [g_0])=\left\{\begin{array}{ll}
\ds\vs \inf _{g\in {\cal C}_2}\tilde { \cal F}_2(g), & \hbox{ if }
\ds {\cal C}_2\neq \emptyset,\\
+\infty,  &\hbox{ if } {\cal C}_2= \emptyset,\end{array}\right.\]
where $\ds \tilde { \cal F}_2(g)=vol(g)^{-\frac{n-4} n}\int_M \s_2(g)\,
dvol(g)$. This is a natural generalization of the Yamabe constant and 
 was considered in \cite{GW3} in the fully 
nonlinear context.

We first prove the following proposition.

\begin{pro}\label{pro1} Let $(M^n,g_0)$ be a compact, 
oriented Riemannian manifold of dimension $n>4$ with $g_0\in \Gamma_2^+$.
The $\s_2$-Yamabe problem is solvable, provided that
\beq\label{0.3}
Y_2(M,[g_0])<Y_2(\S^n).\eeq \end{pro}

The idea to prove the Proposition is  a ``blow-up" analysis, which is a typical 
tool in the field of
semilinear equations. The observation that the fully nonlinear 
Equation (\ref{0.1}) also admits a  blow-up analysis was made in \cite{GW1}.
Inspired by  Yamabe's approach (see  e.g. \cite{Yamabe}, \cite{Aubin1}), 
we first prove the existence of solutions to a ``subcritical" Equation 
(\ref{eq1-e}) for any small $\e>0$.
To prove the existence of solutions of  (\ref{eq1-e}), we use a fully non-linear flow (\ref{flow1}).
We show that this flow globally converges to a solution $u_\e$
 of the subcritical Equation (\ref{eq1-e}).
In fact, $u_\e$ is a minimizer for a corresponding functional. 
This is one of crucial points of this paper.
Then we consider the sequence $u_\e$ as $\e\to 0$. 
Using the blow-up analysis developed in \cite{GW1} and the 
classification of  ``bubbles"  in \cite{CGY3} or \cite{LiLi}, 
we can show that the  sequence
 $u_\e$ subconverges to a solution of (\ref{0.2}) under the
  condition (\ref{0.3}). The flow method to attack the existence of  
  fully nonlinear equations was used by many mathematicians, see for instance
  \cite{Chou}, \cite{Wangxj}, \cite{TW} and \cite{ChouWang}. 
  In the fully nonlinear conformal equations, it was used in \cite{GW2} and
  \cite{GW3}.

Then we show
\begin{pro}\label{pro2} Let $(M^n,g_0)$ be a compact, oriented Riemannian manifold with $g_0\in \Gamma_2^+$.
When $n>8$ and  the Weyl tensor $W_{g_0}\not =0$,
\[Y_2(M,[g_0])<Y_2(\S^n).\] \end{pro}

The proof of this proposition is a delicate gluing argument. We need to construct
suitable test metrics as in \cite{Aubin} and \cite{Schoen} for the ordinary Yamabe problem.
A subtle point in the gluing is that all metrics we constructed should lie in $\Gamma_2^+$.
Recall that in the ordinary Yamabe problem, the test metrics constructed by Aubin and Schoen
have negative scalar curvature somewhere.
To overcome this difficulty, we adopt a method of Gromov-Lawson in their
 construction of metrics of positive scalar curvature. A similar
 method was also used in
 in  \cite{MiWang} for metrics
 of positive isotropic curvature and \cite{GLW} for  metrics of positive $\Gamma_k$-curvatures
 on \lcf manifolds. See also \cite{Haber} and \cite{ShaYang}.
  We believe that by a similar, but  more delicate construction one can prove Proposition 
 \ref{pro2} for $n=8$. For $n=5,6,7$, this problem becomes delicate. 
 We will consider these cases later.
  
  By-products of our work for  flow (\ref{flow1}) are the Poincar\'e type inequality and
 Sobolev   inequality for the operator $\s_2(\uuu)$. (In Section 4, we provide another proof.)

\begin{pro}\label{P0}
Let $(M,g_0)$ be a compact, oriented
 Riemannian manifold with $g_0\in \Gamma_2^+$ and the dimension $n>4$. Then there exists a positive
 constant $\l_1>0$ depending only on $(M,g_0)$ such that for any
 $C^2$ function $u$ with $e^{-2u}g_0\in {\cal C}_2([g_0])$ we have
 \[\int_M \s_2(e^{-2u}g_0)dvol(e^{-2u}g_0) \ge \l_1 \int e^{4u} dvol(e^{-2u}g_0).
 \] 
 Equivalently, for such a function $u$ we have
 \[\int_M e^{(4-n)u}\s_2(\uuu)dvol(g_0) \ge\l_1\int  e^{(4-n)u} dvol(g_0).\]
\end{pro}
 
 \begin{thm}\label{thmS0} Let $(M,g_0)$ be a compact, oriented
 Riemannian manifold with $g_0\in \Gamma_2^+$ and the dimension $n>4$. Then there exists a positive
 constant $C>0$ depending only on $(M,g_0)$ such that for any
 $C^2$ function u with $e^{-2u}g_0\in {\cal C}_2([g_0])$ we have
\[ \int_M \s_2(e^{-2u}g_0)dvol(e^{-2u}g_0) \ge C 
vol(e^{-2u}g_0)^{\frac {n-4} n}.\]
 Equivalently, for such a function $u$ we have
\[ \int_M e^{(4-n)u}\s_2(\uuu)dvol(g_0) \ge C 
(\int_M e^{-nu}dvol(g_0))^{\frac {n-4} n}.\]
\end{thm}

The Sobolev inequality and other geometric inequalities, the Moser-Trudinger
inequality and a conformal quermassintegral inequality
for $\s_k(\uuu)$ for a \lcf manifold were established in \cite{GW3}.
See  also \cite{GVW} and \cite{Gm}.

The method presented here works for a conformal quotient
equation
\[\frac{\s_2(g)}{\s_1(g)}=c,\]
on a general manifold \cite{GeWang2}. See other results for conformal quotient
equations in \cite{GW3}, \cite{GLW2} and \cite{GV2}.

The paper is organized as follows. In Section 2, we discuss various fully nonlinear flows and we prove
local estimates for these flows in Section 3. In Section 4, we establish the Poincar\'e and
Sobolev inequalities. We  prove the global convergence of these fully nonlinear flows and Proposition
\ref{pro1} in Section 5. In Section 6, we prove Proposition \ref{pro2}, and hence 
Theorem \ref{maintheorem}.

\section{Various flows and Ideas of Proof}
Consider the following functional
\begin{equation}\label{eq3}
{ \cal F}_k(g)=\int_M \sk(g)\, dvol(g)\end{equation} and its
normalization $\tilde {\cal F}_k$
\begin{equation}\label{eq3.1}
\tilde { \cal F}_k(g)=vol(g)^{-\frac{n-2k} n}\int_M \sk(g)\,
dvol(g).\end{equation} When $k=2$ or the underlying manifold is
locally conformally flat, Viaclovsky proved that critical points
of $\tilde {\cal F}_2$ are solutions of (\ref{0.1}). Therefore, in these
cases, (\ref{0.1}) is a variational problem. The case when the underlying
manifold is locally conformally flat was studied in \cite{GW2} and \cite{LiLi},
as mentioned in the Introduction.
See also \cite{BV}. 
In this paper we  only consider the case $k=2$. Since the case
$k=2$ and $n\le 4$ was solved in \cite{CGY1}, \cite{GV1} and  \cite{GV2},
we focus on the case $k=2$ and $n>4$.

Recall that  ${\cal C}_2=\{g\in [g_0] \,|\,g\in \Gamma_2^+\}$ and the
Yamabe type constant is defined by
\[Y_2(M, [g_0])=\left\{\begin{array}{ll}
\ds\vs \inf _{g\in {\cal C}_2}\tilde { \cal F}_2(g), & \hbox{ if }
\ds {\cal C}_2\neq \emptyset,\\
\infty,  &\hbox{ if } {\cal C}_2= \emptyset.\end{array}\right.\]
 Our main aim of this
paper is to show that $Y_2(M, [g_0])$ is achieved for  non locally conformally flat manifolds when ${\cal
C}_2([g_0])\neq \emptyset$. In order to achieve our aim, we will first  consider 
 {\it subcritical } equations.
 \beq\label{eq1-e} 
 \s_2^{1/2}(\uuu)=c e^{(\e-2)u},\eeq
 for $\e\in (0, 2]$ and the positive constant $c$.
Its corresponding functional is
\begin{equation}\label{func-e}
\tilde { \cal F}_{2,\e}(g)=V_\e(g)^{-\frac{n-4} {n-2\e}}\int_M \s_2(g)\,
dvol(g),\end{equation} 
where
\[V_\e(g) :=\int _M e^{2\e u} dvol(g)=\int _M e^{(2\e-n) u} dvol(g_0),\]
for $g=e^{-2u}g_0$. It is clear that $V_0(g)=vol(g)$, the volume of $g$ and
$V_2(g)=\int e^{(4-n)u} dvol(g)$.
Set
\[Y_\e(M, [g_0])=  \inf _{g\in {\cal C}_2}\tilde { \cal F}_{2,\e}(g).\]
We will show that $Y_\e(M, [g_0])$ is achieved at $u_\e$, which
is clearly a solution of (\ref{eq1-e}). To prove this we consider 
the following fully nonlinear flow

\begin{equation}\label{flow1}\ba{rcl}\ds\vs
2\frac {du}{dt} &=& \ds-  g^{-1} \cdot \frac {d}{dt}g\\
&=&\ds \left(h(e^{-2u}\s_2^{1/2}(g))- h(r^{1/2}_{\e}(g) e^{(\e-2)u})\right)-s_{\e}(g), \\
&=&\ds  h(\s_{2}^{1/2}(\uuu))-h(r^{1/2}_\e(g)e^{(\e-2)u} )-s_{\e}(g),\ea
\end{equation}
with initial value $u(0)=1$, where $r_\e(g)$ and $s_\e(g)$ are given by for any $\e\in [0,2]$
\[ r_\e(g):=  
\frac {\int_M  \s_2(g)\, dvol(g)}{\int_M  e^{2\e u}
\, dvol(g)}\]
\[ s_\e(g):= \ds \frac {\int_M e^{2\e u}\left(h(e^{-2u} \s^{1/2}_2(g))- 
h(r_\e^{1/2}(g)e^{(\e-2)u}) \right) \, dvol(g)}
{\int_M  e^{2\e u}\, dvol(g)}\] 
and $h:\R_+\to\R$ is smooth  concave function with $h'(t)\ge 1$ for $t\in\R_+$
satisfying
\[h(s)=\begin{cases} 2\log s & \hbox{ if }  t\le 1 \\
s & \hbox{ if }  t\ge 2. \end{cases}\]

  Flow (\ref{flow1}) preserves $V_\e$ and
non-increases ${\cal F}_2$.

\begin{lem}\label{presv}
For any $\e\in [0,2]$, the flow (\ref{flow1}) preserves the functional $V_\e$ and
nonincreases ${\cal F}_2$. In fact,
we have
\begin{equation}
\label{func}
\frac d{dt} {\cal F}_2(g) = -\frac{n-4}{ 2} 
\int_M \left(h(e^{-2u}\s_2^{1/2}(g))- h(r_\e^{1/2}(g)e^{(\e-2)u })\right)( 
\s_2(g)-r_\e
e^{2\e u}) \,dvol(g).
\end{equation}
Moreover, $r_\e$ is bounded.
\end{lem}

\pr We note that 
$$
\frac d{dt}{\cal F}_2(g) = \frac{n-4}{ 2} \int_M (g^{-1} \cdot \frac {d}{dt}g)
\s_2(g) \,dvol(g)
$$
and
$$
\frac d{dt}V_\e(g)= \frac{n-2\e}{ 2} 
\int_M (g^{-1} \cdot \frac {d}{dt}g)e^{2\e u} \,dvol(g)=0.
$$
See the proof in \cite{GW2}.
It is clear that $V_\e$ is preserved along the flow.
 On the other hand, a direct  computation gives
\beq\label{dec}
\begin{array}{ll}
\ds\frac d{dt}{\cal F}_2(g) &=\ds \vs \frac{n-4}{ 2} \int_M (g^{-1} \cdot \frac {d}{dt}g)
(\s_2(g)- r_\e
e^{2\e u}) \,dvol(g)\\
&=\ds -\frac{n-4}{ 2}\ds\int_M \left(h(e^{-2u}\s_2^{1/2}(g))- h(r_\e^{1/2}(g)e^{(\e-2)u })\right)( 
\s_2(g)-r_\e
e^{2\e u}) \,dvol(g),
\end{array}
\eeq
where in the second equality we have used the fact 
$$
\int_M ( \s_2(g)-r_\e
e^{2\e u}) \,dvol(g)=0.
$$
Now it is easy to see that $r_\e$ is bounded. 
\qed



In fact, flow (\ref{flow1}) strictly decreases the functional
${\cal F}_2$ except at the solutions of  Equation (\ref{eq1-e}).
When $\e=0$ Equation (\ref{eq1-e}) is just
(\ref{0.2}). When $\e=2$ Equation (\ref{eq1-e}) is a corresponding
 equation for a nonlinear eigenvalue problem, which was considered
 in \cite{GW2.1}. See also  Section 4.

Since $g_0 \in \Gamma_2^+$, flow  (\ref{flow1}) is parabolic near
$t=0$. By the standard implicit function theorem we have
the following short-time existence result.
\begin{pro}\label{pro4} For any $g_0\in C^2(M)$ with $g_0 \in
\Gamma^+_2$, there exists a positive constant $T^*\in (0, \infty]$ such that
 flow (\ref{flow1}) exists and is parabolic for $t \in [0, T^*)$,
and $\forall T<T^*,$
\[ g \in C^{3,\alpha}([0,T] \times M), \forall 0<\alpha <1,
 \quad \text{and} \quad g(t)\in \Gamma^+_2 .\]
 \end{pro}

We assume that $T^*$ is the largest number, for which Proposition \ref{pro4}
holds.
 We first show that the global convergence of flow (\ref{flow1}) when
 $\e=2$. The global convergence implies a Poincar\'e type
 inequality. Then, using this inequality and the divergence free of
 the first Newton transformation of the Schouten tensor, which was
 an observation  Viaclovsky, we obtain an optimal Sobolev
 inequality. By establishing a flow version of local gradient
 estimates, which was proved in \cite{GW1},  we show that flow (\ref{flow1}) 
 globally
 converges to a solution $u_\e$ of (\ref{eq1-e}) for any $\e \in
 (0,2]$. With the help of the local estimate obtained in
 \cite{GW1} and a classification in \cite{LiLi} or \cite{CGY3}, we
 show that $u_\e$ subconverges to a solution $u_0$ of (\ref{0.2}),
 provided that
 \beq \label{Aubin} Y_2(M,[g_0])<Y_2(\S^n).\eeq
 In this case, it is clear that $u_0$ is the minimum of
 $\tilde{\cal F}_2$. 

\section{Local Estimates}
In this section, we will establish a local estimate for solutions
of (\ref{flow1}), which is a parabolic version of 
a local estimate for solutions
of (\ref{0.2}) obtained in \cite{GW1}.

\begin{thm}\label{thmlocal1} Let $u$ be a solution of (\ref{flow1})
with $\e\in [0,2]$ in a geodesic ball $   B_r\times[0, T]$ for
$T<T^*$ and $r<r_0$, the injectivity radius of $M$. There is a
constant $C>0$ depending only on  $(B_r, g_0)$ such that
for any $(x,t) \in B_{r/2}\times[0, T]$
\begin{equation}\label{eq_thm1.1}
|\n u|^2+|\n^2 u|\le C(1+  e^{-(2-\e)\inf_{(x,t)\in B_r\times[0, T] }
u(x,t)}).
\end{equation}
\end{thm}\

\pr 
 The proof follows \cite{GW1} closely. We only point out the different places. Without loss of generality, we assume $r=1$.
Let $ \rho\in C^\infty_0(B_1)$ be a test function defined as in
\cite{GW1}. such that
\begin{equation}\begin{array}{rcll}\label{8}
\vs \ds \rho & \ge & 0, & \hbox{ in } B_1,\\
\vs\ds \rho & =& 1, &\hbox{ in } B_{1/2},\\
\vs\ds |\n \rho (x)|& \le &2 b_0 \rho^{1/2}(x), & \hbox{ in } B_1,\\
|\n^2 \rho| & \le  & b_0,   & \hbox{ in } B_1.
\end{array}\end{equation}
Here $b_0>1$  is a constant. Set $H(x, t)=\rho |\n u|^2$. Let
$(x_0, t_0)$ be the maximum of $H$ in $M\times [0,T]$. Without loss of generality, we assume $t_0>0$. We have at
$(x_0, t_0)$ that
\begin{eqnarray}
\label{a11}  \ds\vs 0 & \le & H_t=\ds 2\rho \sum_{l} u_l u_{lt} ,\\
 \label{a12} 0& = & \ds\vs H_j= \rho_j|\n u|^2+ 2\rho \sum_{l} u_l u_{lj}, \\
 \label{a13} 0 & \ge & (H_{ij}).
\end{eqnarray}

 Let $W=(w_{ij})$ be an $n\times n$ matrix
 with $w_{ij}=\n ^2_{ij} u+u_iu_j-\frac {|\n u|^2}2
 (g_0)_{ij}+(S_{g_0})_{ij}$. Here $u_i$ and $u_{ij}$ are the first
 and second derivatives of $u$ with respect to the background metric $g_0$.
By choosing suitable normal coordinates, we may assume  that $W$ is
diagonal at $(x_0,t_0)$, and hence we have at $(x_0, t_0)$,
\begin{equation}\begin{array}{rcll}\label{9}
\vs \ds w_{ii}&=&\ds u_{ii}+u_i^2-\frac 12 |\n u|^2+(S_{g_0})_{ii},\\
u_{ij}&=&-u_iu_j-(S_{g_0})_{ij}, \quad \quad \forall i \neq j.
\end{array}\end{equation}
In view of (\ref{8}), (\ref{a12}) and (\ref{9}), we have  at
$(x_0, t_0)$
 \begin{equation} \label{11}
 |\sum_{l=1}^n u_{il}u_l|\le  b_0 \rho^{-1/2}|\n u|^2.
 \end{equation}
 We may assume that
\[H(x_0, t_0)\ge A^2_0b^2_0,\]
 i. e.,  $\rho^{-1/2}\le  \frac 1{ A_0b_0}|\n u| $, and
\[|\nabla S_{g_0}|+ |S_{g_0}|\le A_0^{-1}|\n u|^2,\]
where $A_0>1$ is a large, but fixed number to be chosen later,
otherwise we are done.
  Thus, from (\ref{11}) we have
\begin{equation}
\label{12}|\sum_{l=1}^n u_{il}u_l| \le \frac {|\n u|^3}
{A_0}(x_0, t_0).\end{equation} Set $F=h(\s^{1/2}_2(W))$ and
\[F^{ij}=\frac {\partial  F(W)}{\partial w_{ij}}.\]
Note that flow (\ref{flow1}) is equivalent to
 $2u_t=F-h(r_\e^{1/2}e^{(\e-2)u})-s_\e^{1/2}$ and
 $F^{ij}$ is diagonal at $(x_0, t_0)$.
Since matrix $(F^{ij})$ is positive definite, {from} (\ref{a11})
and (\ref{a13}) we have
\begin{equation}\label{rcl}\ba{rcl}
0& \ge & \ds\vs \sum_{i,j}F^{ij}H_{ij}-2H_t\\
&=& \ds\vs \sum_{i,j} F^{ij}
\left\{\left(-2\frac{\rho_i\rho_j}\rho+\rho_{ij}\right)|\n u|^2
+2\rho \sum_{l}u_{lij}u_l+2\rho \sum_{l} u_{il}u_{jl}\right\}-4\rho \sum_{l} u_lu_{lt}.
\ea\eeq

 We need to estimate
 the term $\ds\sum_{i,j,l}F^{ij}
 u_{lij}u_l-2 \sum_{l} u_lu_{lt}$. Since  changing the order of
 derivatives only causes a low order term, we have
 \beq\label{u1}
 \ba{rcl}
 \ds\vs\sum_{i,j,l} F^{ij}
 u_{lij}u_l-2\sum_{l} u_lu_{lt}&\ge
 & \ds \sum_{i,j,l} F^{ij}
 u_{ijl}u_l-2\sum_{l}u_lu_{lt}-c\sum_i F^{ii} |\n u|^2\\
 \vs &\ge & \ds\sum_{i,j,l} F^{ij} (w_{ij})_lu_l
 -\sum_{i,l}F^{ii}(u^2_i-\frac 12 |\n u|^2)_l u_l-2\sum_{l}u_lu_{lt}\\
 \vs&&\ds-c\sum_i F^{ii} |\n u|^2-\sum_{i,l}F^{ii}\nabla_l(S_{g_0})_{ii}u_l \\
 \ds\vs &\ge & \ds \sum_l(F_l-2u_{tl})u_l-c{A_0}^{-1}\sum_i F^{ii} |\n u|^4\\
\ds\vs &\ge &\ds (\e-2) h'(r_\e^{1/2}e^{(\e-2)u})r^{1/2}_\e e^{(\e-2)u} 
 |\n u|^2-c{A_0}^{-1}\sum_i F^{ii} |\n
u|^4, \ea\eeq where we have used (\ref{flow1})  and (\ref{12}).
Here  $c$ is a constant independent of $u$, but it may vary from
line to line. The term
 $\ds\sum_{i,j}F^{ij}
(-2\frac{\rho_i\rho_j}\rho+\rho_{ij})|\n u|^2$ is  bounded from below by
$\ds -10 b_0^2 |\n u|^2\sum_{j}F^{jj}$. For the term
$F^{ij}u_{il}u_{jl}$ we have the following  crucial Lemma.
\begin{lem}[\cite{GW1}]\label{lem1} There is a constant $A_0$ sufficient
large (depending only on $n$, and $\|g_0\|_{C^3(B_1)}$), such
that,
\begin{equation}\label{15}\sum_{i,j,l}F^{ij}u_{il}u_{jl}
\ge A_0^{-{\frac 34}}|\n u|^4\sum_{i\ge 1}  F^{ii}.
\end{equation}
\end{lem}
Altogether gives us \beq\label{z1} (A_0^{-\frac 34}-c
A_0^{-1})\rho |\n u|^4 \sum_i F^{ii} \le 10 b_0^2|\n u|^2 \sum_i
F^{ii}+c\rho(1+ e^{(\e-2)u} )
 |\n u|^2.\eeq 
 By the Newton-McLaurin inequality and the fact that $h'(t)\geq 1$ for any $ t\geq 0$,  it is easy
to check that
 \[\sum_i F^{ii}\ge 1,\]
 which, together with (\ref{z1}), proves the local gradient estimate
 \[ |\n u|^2 \le C(1+ e^{(2-\e)\inf_{(x,t)\in B_r\times [0,T]}u(x,t)}),\]
 for some constant $C>0$ depending only on $(B_r, g_0)$.
 
 \
 
 Now we show the local estimates for second order derivatives.
 Since $e^{-2u}g_0\in \Gamma_2^+$, to bound $|\n^2 u|$ we only need bound
 $\D u$ from above. This is a well-known fact, see for instance \cite{GW1}.
 Set 
 \[G=\rho(\D u+|\n u|^2),\]
 where $\rho$ is defined as above. Let $(y_0,t_0)$ be a maximum point
 of $G$ in $M\times [0,T]$. 
 Without loss of generality, we assume $G(y_0, t_0)>1+2\max H(x,t)$, 
 $t_0>0$ and $(u_{ij})$ is diagonal at $(y_0,t_0)$. Recall that 
 $H=\rho |\n u|^2$. Hence we have
 \[ 0<\rho \D u(y_0,t_0) \le G(y_0,t_0) \le 2\rho \D u(y_0,t_0	).\]
At $(y_0, t_0)$, we have
\begin{eqnarray}\label{y1}
0&\le & G_t=\rho\sum_{l}(u_{llt}+2u_lu_{lt}),\\
\label{y2} 0&=&G_j=\frac{ \rho_j }{\rho} G + \rho\sum_{l\ge 1}
 (u_{llj}+2u_l u_{lj}),
\quad \hbox{ for any } j,\\
0&\ge &G_{ij} =\ds \frac{ \rho \rho_{ij}-2\rho_i\rho_j}{\rho^2} G+
\rho \sum_{l\ge 1} (u_{llij}+2u_{li}u_{lj}+2 u_l u_{lij}).\end{eqnarray} 
Recall that $F^{ij}=\frac{\partial}{\partial w_{ij}}F$  
is non-negative definite. Hence, we have
\[\ba{rcl}
0& \ge &\ds\vs  \sum _{i,j\ge 1}F^{ij}G_{ij}
-2G_t\\
& \ge &  \ds\vs\sum_{i,j\ge 1}
 F^{ij} \frac{ \rho \rho_{ij}-2\rho_i\rho_j}{\rho^2} G +
\rho \sum_{i,j, l\ge 1} F^{ij}(u_{ijll}+2u_{li}u_{lj}+2 u_l u_{lij})\\
&& \ds-2\rho\sum_{l}(u_{llt}+2u_lu_{lt}) -C\rho\sum_i(|u_{ii}|+|u_{i}|)\sum_{i,j}|F^{ij}|,\ea \]
where the last term comes from the commutators related to
the curvature tensor of $g_0$ and its derivatives.
 From the definition of
 $\rho$, we have
\[ \sum_{i,j\ge 1}
 F^{ij} \frac{ \rho \rho_{ij}-2\rho_i\rho_j}{\rho^2} G \ge -C
 \sum_{i,j \ge 1}|F^{ij}|\frac 1 \rho  G.\]
By the concavity of $\s_2^{1/2}$, we have
\beq\label{x1}\ba{rcl} 
\ds\vs  \sum_{i,j, l\ge 1} F^{ij} u_{ijll} &=& \ds    
\sum_{i,j, l\ge 1} F^{ij} w_{ijll} -
\sum_{i,j, l\ge 1} F^{ij} (u_iu_j-\frac 12 |\nabla u|^2 (g_0)_{ij}+(S_{g_0})_{ij})_{ll}\\
&\ge & \ds\vs \sum_l F_{ll}-
\sum_{i,j, l\ge 1} F^{ij} (u_iu_j-\frac 12 |\nabla u|^2 (g_0)_{ij}+(S_{g_0})_{ij})_{ll}.
\ea\eeq
We also have
\beq\label{x2}\ba{rcl}
\ds\vs \sum_{i,j, l} F^{ij} u_{l}u_{lij} &=&\ds
 \sum_{i,j, l} F^{ij} u_{l}(w_{ij})_l
 -\sum_{i,j, l} F^{ij}u_l(u_iu_j-\frac 12 |\nabla u|^2 (g_0)_{ij}+(S_{g_0})_{ij})_{l}
 \\
 &&\ds +\sum_{i,j, l} F^{ij}u_l(u_{lij}-u_{ijl})\\
 &=& \ds \sum_{l}F_lu_l-\sum F^{ij}u_l(u_iu_j-\frac 12 |\nabla u|^2 (g_0)_{ij}+(S_{g_0})_{ij})_{l}\\
 &&\ds +\sum_{i,j, l} F^{ij}u_l(u_{lij}-u_{ijl}).\ea\eeq
Hence, we have
\beq\label{x3}\ba{rcl} 
\ds\vs \sum_{i,j, l\ge 1} F^{ij}(u_{ijll}+2u_{li}u_{lj}
+2 u_l u_{lij})
&\ge &\ds\vs \sum_l (F_{ll}+2F_lu_l)-2\sum_{i,j,l} F^{ij}u_iu_{jll}+\sum_{j,k,l} 
F^{jj}u_ku_{kll}\\
&&\ds\vs -2\sum_{i,j,l} F^{ij}u_l(u_iu_j-\frac12 |\nabla u|^2 (g_0)_{ij}+(S_{g_0})_{ij})_l\\
&&\ds +\sum_{i,k,l}F^{ii}(u_{kl} )^2-C(1+\frac G\rho)\sum_{i,j} |F^{ij}|.
\ea\eeq
{From} (\ref{y1}) and Equation (\ref{flow1}), we have
\beq\label{x4}\ba{rcl}
\ds\vs  \rho\sum_{l}(F_{ll}+2F_lu_l)&\ge &\ds
2\rho \sum _l(u_{tll}+2u_lu_{tl})-C (2-\e)G(1+ e^{(\e-2)u}).\ea
\eeq
The term $\ds -2\sum_{i,j,l} F^{ij}u_iu_{jll}+\sum_{j,k,l} 
F^{jj}u_ku_{kll}$ can be controlled as in \cite{GW1} with the help of
(\ref{y2}). And the other terms in (\ref{x3}) can easily be estimated. On the other hand, it follows from the positivity of $(F^{ij})$ that
\[
\sum_{i,j} |F^{ij}|\leq C \sum_{i} F^{ii}
\]
This completes the proof of the Theorem.\qed

 {From} the local estimates, we have
\begin{cor}\label{coro1} If ``bubble" occurs, i.e., $\inf_{M\times [0,T^*)} u=
 -\infty$, then
there is a positive constant $c_0>0$ such that
\[\lim_{\d\to 0}\lim_{t\to T^*}
V_\e(g, B_\d)>c_0.\]
\end{cor}


\section{A Poincar\'e inequality and a Sobolev inequality}

The Sobolev inequality is a very important analytic tool in many
problems arising from  analysis and geometry. It plays a crucial
role in the resolution of the Yamabe problem, which was solved
completely by Yamabe \cite{Yamabe}, Trudinger \cite{Trudinger},
Aubin \cite{Aubin} and Schoen \cite{Schoen}. 
See various optimal
Sobolev inequalities in \cite{Hebey}.
 In this section we are
interested in a similar type inequality for the class of a fully
nonlinear conformal  operators $\s_k(\uuu)$. In \cite{GW3},
the Sobolev inequality was generalized to the operator
$\s_k(\uuu)$ for $k<n/2,$ if the underlying manifold is locally
conformally flat. Namely,
 \begin{thm}[\cite{GW1}]\label{thmS} Let $(M^n,g_0)$ be a compact, oriented
 Riemannian manifold with $g_0\in \Gamma_k^+$ and $k<n/2$.
 Assume that $(M, g_0)$ is locally conformally flat,
 then there exists a positive
 constant $C>0$ depending only on $n$, $k$ and $(M,g_0)$ such that for any
 $C^2$ function u with $e^{-2u}g_0\in {\cal C}_k([g_0])$ we have
 \begin{equation}\label{S00}
 \int_M \s_k(e^{-2u}g_0)dvol(e^{-2u}g_0) \ge C vol(e^{-2u}g_0)^{\frac {n-2k} n}.
 \end{equation}Equivalently, for such a function $u$ we have
  \begin{equation}\label{S002}
 \int_M e^{(2k-n)u}\s_k(\uuu)dvol(g_0) \ge C (\int_M e^{-nu}dvol(g_0))^{\frac {n-2k} n}.
 \end{equation}
\end{thm}
When $k=1$, inequality (\ref{S}) is just the Sobolev inequality.
The proof of Theorem \ref{thmS} uses a Yamabe type flow.
See also the work of \cite{Gm}.

In this section, we  establish the Sobolev inequality for $k=2$
without the flatness condition.
 \begin{thm}\label{thmS2} Let $(M,g_0)$ be a compact, oriented
 Riemannian manifold with $g_0\in \Gamma^+_2$ and the dimension $n>4$. Then there exists a positive
 constant $C>0$ depending only on $(M,g_0)$ such that for any
 $C^2$ function u with $e^{-2u}g_0\in {\cal C}_2([g_0])$ we have
 \begin{equation}\label{S}
 \int_M \s_2(e^{-2u}g_0)dvol(e^{-2u}g_0) \ge C vol(e^{-2u}g_0)^{\frac {n-4} n}.
 \end{equation}Equivalently, for such a function $u$ we have
  \begin{equation}\label{S2}
 \int_M e^{(4-n)u}\s_2(\uuu)dvol(g_0) \ge C (\int_M e^{-nu}dvol(g_0))^{\frac {n-4} n}.
 \end{equation}
\end{thm}

First we prove a Poincar\'e type inequality, which will be used in
the proof of our Sobolev inequality. The usual Poincar\'e type
inequality is associated to the first eigenvalue problem. In our
case, there is a nonlinear eigenvalue problem, which was studied
in \cite{GW2.1}.
\begin{pro}\label{pro5} Let $(M,g_0)$ be a compact manifold with
$g_0\in \Gamma_k^+$. Then there is a function $u$ with
$e^{-2u}g_0\in \Gamma_k^+$ satisfying \beq\label{eigenvalue}\s_k
(\uuu)=\l_1>0.\eeq Moreover the constant $\l_1$ is unique and the solution is unique up to a
constant.
\end{pro}
An elliptic method was used in the proof, which was motivated by
a method introduced in \cite{PL-Lions}. See also \cite{Wangxj} for
a Hessian operator. In view of Proposition \ref{pro5}, one may
guess that
\begin{equation}\label{2.1}\int e^{2k u} \s_k(\uuu) dvol(e^{-2u}g_0)\ge \l_1 \int
e^{2ku} dvol(e^{-2u }g_0),\end{equation}
 for any $u$ with
$e^{-2u}g_0\in\Gamma_k^+$. It is easy to see that when $k=1$
inequality (\ref{2.1}) holds. In fact it is the Poincar\'e
inequality. In this section, we show that (\ref{2.1}) holds for
$k=2$ by flow (\ref{flow1}) with $\e=2$.

\begin{pro}\label{P}
Let $(M,g_0)$ be a compact, oriented
 Riemannian manifold with $g_0\in \Gamma_2^+$ and the dimension $n>4$. Then for any
 $C^2$ function $u$ with $e^{-2u}g_0\in {\cal C}_2([g_0])$ we have
 \begin{equation}\label{P1}
 \int_M \s_2(e^{-2u}g_0)dvol(e^{-2u}g_0) \ge \l_1 \int e^{4u} dvol(e^{-2u}g_0).
 \end{equation}Equivalently, for such a function $u$ we have
  \begin{equation}\label{P2}
 \int_M e^{(4-n)u}\s_2(\uuu)dvol(g_0) \ge\l_1\int e^{(4-n)u} dvol(g_0).
 \end{equation}
\end{pro}

\pr To prove the Proposition, we consider flow  (\ref{flow1})
with $\e=2$. 
We want to show that the flow converges
globally to a solution obtained in Proposition \ref{pro5}.
 By Theorem \ref{thmlocal1}, we have
 \beq\label{a1}
 |\n^2 u|+|\n u|^2(x,t)\le C,\eeq
 where $C$ is a constant independent of $(x,t)\in M\times
 [0,T^*)$. Since the flow preserves the functional $V_2$, in view of (\ref{a1})
  we have that $|u|\le C$, for some constant $C>0$. 
  Now following the method in \cite{GW2}
  we can show that
  \[\s_2(g)>c_0,\]
  for some constant $c_0$ independent of $t$. See the proof in the next section.
   Hence, this flow exists globally and is
  uniformly elliptic.
  By the result of Krylov, $g(t)\in C^{4+\a, 2+\a}$. Since the
  flow satisfies (\ref{func}), one can show that for any sequence of
  $\{t_i\}$ with $t_i\to\infty$ there is a subsequence, still denoted by
  $\{t_i\}$, such that $g(t_i)$ converges strongly to $g^*$, which
  satisfies (\ref{eigenvalue}). On the other hand, $V_2(g^*)\equiv V_2 (g(t))$. By the uniqueness in Proposition \ref{pro5}, 
   one can  show that the flow globally
  converges to $g^*$.
  Since the flow preserves $V_2$ and decreases ${\cal F}_2$, we
  have
  \[{\cal F}_2(g) \ge {\cal F}_2(g^*),\]
  for any $g\in {\cal C}_2$. This is the Poincar\'e inequality that
  we want to prove.
  \qed
\

\noindent{\it Proof of Theorem \ref{thmS2}.} Let $g=e^{-2u}g_0$. We
have
\[
2\s_2=\sum_{i,j}T^{ij}S_{ij},
\]
where $T(g)^{ij}= \s_1(g)g^{ij}-S(g)^{ij}$ is the so-called the first
Newton transformation. We will use the following formulas \beq
S(g)_{ij}=u_{ij}+u_iu_j-\frac 12|\n u|^2_{g_0} (g_0)_{ij}+S(g_0)_{ij}
\eeq and \beq \widetilde \n^2_{ij}u =u_{ij}+2u_iu_j-|\n
u|^2_{g_0} (g_0)_{ij}, \eeq where $\widetilde \n$ are the
derivatives w. r. t $g$. Thus, \beq 2\s_2(g)=\sum_{i,j}T(g)^{ij}\widetilde
\n^2_{ij} u- \sum_{i,j}T(g)^{ij}u_iu_j+\frac {n-1}2 \s_1(g)|\widetilde\n
u|^2_g+\sum_{i,j}T(g)^{ij}S(g_0)_{ij}. \eeq Here we have used $\tr
T(g)=(n-1)\s_1(g)$. Note that \beq \label{positivity} \sum_{i,j}T(g)^{ij}
S(g_0)_{ij}>0, \eeq thanks to Garding's inequality \[
\sum_{i,j}T(g)^{ij}S(g_0)_{ij}\ge 2e^{2u}\s_2^{1/2}(g)\s_2^{1/2}(g_0).\] Due to an
observation of Viaclovsky, $\sum_{i}\widetilde \n_i T(g)^{ij}=0$,
we have 
\beq \ba{rcl}\label{45}
\ds\vs 2\int \s_2(g) &=&  \ds-\int \sum_{i,j}T(g)^{ij} u_i u_j dvol(g)+\frac {n-1}
2\int \s_1(g)|\widetilde\n u|^2_g dvol(g)\\
&&\ds +\int\sum_{i,j} T(g)^{ij} S(g_0)_{ij}
d vol(g). \ea\eeq Recall that $T(g)=\s_1(g) g-S(g)$. We have \beq
\begin{array}{rcl}
\ds\vs -\int \sum_{i,j}T(g)^{ij}u_iu_j dvol(g)&=&\ds -\int
\s_1(g)|\widetilde\n u|^2_g
dvol(g)+\int \sum_{i,j}S(g)^{ij}u_iu_jdvol(g)\\
&=&\ds\vs -\int \s_1(g)|\widetilde\n u|^2_g dvol(g)+\int
\sum_{i,j}\widetilde \n^{ij} u u_i u_j dol(g)\\ &&\ds\vs -\frac 12 \int
|\widetilde\n u|^4_g dvol(g)+\int\sum_{i,j} S(g_0)^{ij} u_i u_j
dvol(g).\end{array} \eeq and \beq
\begin{array}{rcl}
\ds\vs \int \sum_{i,j}\widetilde \n^{ij} u u_iu_j dvol(g)&=&\ds \frac 12
\int
\sum_{i}\widetilde \n^i(|\widetilde\n u|^2_g) u_idvol(g)\\
&=& \ds\vs -\frac 12\int |\widetilde\n u|^2_g \tr (\widetilde \n^2 u) dvol(g)\\
&=&\ds\vs -\frac 12 \int \s_1(g)|\widetilde\n u|^2_g +\frac {n-2}
4 \int |\widetilde\n u|^4_g
\\
&& +\ds \int\frac 12 \s_1(g_0) |\widetilde\n
u|^2_g e^{2u}dvol(g)\end{array} \eeq 
Hence \beq \label{47}
\begin{array}{rcl}\ds\vs  -\int \sum_{i,j}T(g)^{ij}u_iu_j dvol(g)&=&\ds
-\frac 32 \int |\widetilde\n u|^2_g +\frac {n-4} 4 \int |\widetilde\n u|^4_g \\
&&\ds\vs +\int\sum_{i,j} S(g_0)^{ij}u_iu_j+\frac 12 \int \s_1(g_0)
|\widetilde\n u|^2_g e^{2u},\end{array} \eeq where all integrals
are w.r.t $g$. (\ref{45}) and (\ref{47}) give us \beq
\begin{array}{rcl}
\ds\vs 2\int \s_2(g) dvol(g)&=& \ds \frac {n-4}2 \int \s_1(g)
|\widetilde\n u|^2_g
dvol{g}+\frac {n-4} 4 \int|\widetilde\n u|^4_g dvol(g)\\
&&\ds\vs +\int \sum_{i,j}T^{ij} S(g_0)_{ij} dvol(g)+\int \sum_{i,j}S(g_0)^{ij}u_iu_j
dvol(g)\\
&&\ds+\frac 12 \int\s_1(g_0)|\widetilde\n u|^2_g e^{2u}
dvol(g)\end{array}. \eeq
Finally, we obtain \beq
\begin{array}{rcl}
\ds\vs 2\int \s_2(g) dvol(g)&=& \ds \frac {n-4}2 \int \s_1(g)|\n
u|^2_{g_0}e^{2u}
dvol{g}+\frac {n-4} 4 \int|\n u|^4_{g_0} e^{4u}dvol(g)\\
&&\ds \vs +\int \sum_{i,j}T^{ij} S(g_0)_{ij} dvol(g)+\int \sum_{i,j}S(g_0)^{ij}u_iu_j
dvol(g)\\
&&\ds+\frac 12 \int\s_1(g_0)|\n u|^2_{g_0} e^{4u}
dvol(g).\end{array} \eeq Recall (\ref{positivity}) and positivity
of $\s_1(g)$ and $\s_1(g_0)$. Using the estimates \beq
\sum_{i,j} S(g_0)^{ij}u_iu_j \geq -c |\n u|^2_{g_0}e^{4u}\geq -\frac {n-4} 8
|\n u|^4_{g_0} e^{4u}- \frac {2c^2}{n-4}  e^{4u} \eeq we deduce
\beq \label{52}
\begin{array}{rcl}
\ds\vs 2\int \s_2(g) dvol(g)&\geq& \ds \frac {n-4} 8 \int|\n
u|^4_{g_0} e^{4u}dvol(g) -c \int e^{4u} dvol(g).
\end{array}
\eeq In view of the  Poincar\'e inequality (\ref{P1}), 
 the Sobolev inequality (\ref{S}) follows from (\ref{52}).
 \qed
 
 We remark that a similar method was used to obtain Sobolev inequalities
 on \lcf manifolds by Gonz\'alez in \cite{Gm}. The arguement given in the next section will provide
 another proof of Theorem \ref{thmS2}.

\section{Global convergence of flow (\ref{flow1}) when $\e>0$}

\begin{pro}\label{prop-e} For any $\e\in (0,2]$, flow (\ref{flow1})
converges globally to $u_\e$, which satisfies (\ref{eq1-e}).
\end{pro}

\pr For any $t\in [0,T^*)$, set \[m(t)=\min_{(x,s)\in M\times [0,t]} u(x,s).\]
If $\inf_{t\in [0, T^*)} m(t) >-\infty$, then by estimates given in Section 3,
we have a uniform bound of $|\n u|^2+|\n^2 u|$. Since  flow (\ref{flow1})
preserves the  functional $V_\e$, we have a uniform  $C^2$ bound.
Now we claim that there is a constant $c>0$ such that
\beq\label{claim} F(x,t)\ge c>0, \quad \hbox{ for any } (x,t)\in M\times [0,T^*).
\eeq
Recall that $F=\s_2^{1/2}(\uuu)$. We will prove the 
claim at the end of the proof. (\ref{claim}) implies that
flow (\ref{flow1}) is uniformly elliptic in $M\times [0,T^*)$. Hence, by 
Krylov's result, $u$ has a uniform bound for higher order derivatives, which
implies first that $T^*=\infty$, the global existence.
The global convergence of (\ref{flow1}) with $\e\in (0,2]$
follows now closely the argument presented in \cite{GW2}, which,
in turn, follows closely the argument given in \cite{Simon} and \cite{Andrews}. Therefore, to prove the Proposition, we only need to exclude  that
\beq\label{z03}\inf_{t\in [0, T^*) }m(t) =-\infty.\eeq
We assume by contradiction that $\inf_{t\in [0, T^*)} m(t) =-\infty$.
Let $T_i$ be a sequence tending to $T^*$ with $m(T_i)\to -\infty$ as 
$i\to \infty$. Let $(x_i,t_i)\in M\times [0,T_i]$ with
$u(x_i, t_i)=m(T_i)$. Fixing $\d\in (\frac 25, \frac12)$, we consider 
$r_i=\frac{\e}{2}|m(T_i)|e^{(1-\d\e)m(T_i)}$. Clearly, we have $r_i\to 0$. 
It follows from Theorem \ref{thmlocal1} that for sufficiently large $i$
$$\ba{rcl}
\vs u(x,t_i) &\le & m(T_i)+|\n u|r_i \\
&\le & \vs\ds  m(T_i)+Ce^{(\frac {\e}2-1)m(T_i)}\frac {\e} 2|m(T_i)| 
e^{(1-\d \e)m(T_i)}\\
&= & \vs\ds  m(T_i)+C\frac {\e} 2|m(T_i)| e^{\e(\frac 12 -\d) m(T_i)}\\
&\leq &\ds (1-\kappa)m(T_i),\qquad \forall x\in B(x_i,r_i), \ea
$$
for some $\kappa\in (0,(\d-{\frac 2 n})\e)$. Note that $\d-\frac 2 n >0$, for
$n\ge 5$.
 Therefore, we obtain
\[ 
\begin{array}{lll}
\ds{\int_{B(x_i,r_i)} e^{2\e u}dvol(g)}&\geq \ds{\int_{B(x_i,r_i)} e^{(2\e-n)m(T_i)(1-\kappa)} 
dvol(g_0)\geq Ce^{(2\e-n)m(T_i)(1-\kappa)}r_i^n}\\
&\geq C\left(\frac{|m(T_i)|\e}{2}\right)^n\to \infty.
\end{array}
\]
Hence, this fact contradicts the boundedness of $V_\e$.

Now we remain to prove Claim (\ref{claim}). For any $0<T<T^*$, set
\[
T_1:=\inf\{T'\in [0,T]|\quad \forall(x, t)\in M\times [T',T],\quad r_\e^{1/2}(g(t))e^{(\e-2)u(x,t)}<1/2 \}
\]
It is clear that  $\forall t\in [0,T_1]$ we have $r_\e(g(t))>C$ for positive constant $C>0$ independent of $T$. Let us consider a function $H:M\times [0,T_1]$ defined by 
\[\ba{rrl}
H&:=&\ds\vs \frac12 (h(\s_2^{1/2}(\uuu)-h(r_\e^{1/2}(g)e^{(\e-2)u}))-e^{-u}\\
&=&\ds u_t+\frac12 s_\e(g)-e^{-u}. \ea\]
We first compute the evolution equation for $\s_2$.
A direct computation, see for instance Lemma 2 in \cite{GW2},
gives
\[\ba{rcl}
\ds\vs \frac d{dt}\s_2&=&\ds 
2\s_2g\cdot \frac d{dt}(g^{-1})+\tr\{T_1(S_g) g^{-1}\frac d{dt}S_g\}\\
&=& \ds\vs 4\s_2(g) u_t+\tr\{T_1(S_g) g^{-1} \tilde\n_g^2(u_t)\},\\
\ea\]
where $\tilde \n$ is the derivatives with respect to the evolved metric $g$.
Without loss of generality, we assume that the minimum of $H$ is achieved
at $(x_0, t_0)\in M\times (0, T_1]$. We will show that there is a
constant $c_0>0$ independent of $T$ such that
\beq \label{54.1}\s_2(g)(x_0, t_0) > c_0.\eeq
Since $|u|$ has a uniform bound, without loss of generality
we may assume that at $(x_0,t_0)$
\[e^{-2 u}\s_2^{1/2}(g) <1/2.\] 
Recall that $h(t)=2\log t$ for $t<1$. Hence, in a small neighborhood of
$(x_0,t_0)$
\[H=\log (e^{-2u}\s_2^{1/2}(g))-\frac 12 h(r_\e^{1/2}(g)e^{(\e-2)u})-e^{-u}.\]
 Let us use $O(1)$ denote
terms with a uniform  bound. We have near $(x_0,t_0)$
\beq\label{v1}\ba{rcl}
\ds\vs \frac{d}{dt} H &=&
\ds \frac1 {2\s_2(g)}\tr\{T_1(S_g) g^{-1} \tilde\n_g^2(u_t)\} 
-\frac 14 h'(r_\e^{1/2}(g)e^{(\e-2)u}) \frac {dr_\e(g)} {dt}
r_\e^{-1/2}(g)e^{(\e-2)u}\\
&&\ds\vs +\left[e^{-u}+\frac{2-\e}2 h'(r_\e^{1/2}(g)e^{(\e-2)u})
r_\e^{1/2}(g)e^{(\e-2)u}\right]u_t
\\
&\ge
&\ds\vs \frac1 {2\s_2(g)}\tr\{T_1(S_g) g^{-1} \tilde \n_g^2(H)\}+ 
\frac1 {2\s_2(g)}\tr\{T_1(S_g) g^{-1} 
\tilde\n_g^2(e^{-u})\}\\
&&+ \ds \left[e^{-u}+\frac{2-\e}2 h'(r_\e^{1/2}(g)e^{(\e-2)u})
r_\e^{1/2}(g)e^{(\e-2)u}\right]u_t, \ea
\eeq
where we have used $\frac {dr_\e(g)}{dt}\le 0$. See Lemma \ref{presv}.
Let $F=\log \s_2(\uuu)$ and $F^{ij}=\frac{\partial}{\partial w_{ij}}F$. 
Since $(x_0,t_0)$ is the minimum of $H$ in $ M\times [0, T_1]$,
 at this point,  we  have
\[ \frac {dH}{dt} \le 0, \]
\[0= H_l = \frac 1 2\sum_{ij}F^{ij} w_{ijl}+\left(e^{-u}+\frac{2-\e}2 h'(r_\e^{1/2}(g)e^{(\e-2)u})
r_\e^{1/2}(g)e^{(\e-2)u}\right)u_l=0 \quad \forall l\]
and
\[ (H_{ij}) \text{ is non-negative definite.}\]
Note that 
\[ (\tilde \n^2_g)_{ij} H=H_{ij}+ u_iH_j+u_jH_i-\sum_l u_lH_l \d_{ij}=H_{ij},
\]
at $(x_0,t_0)$, 
where  $H_{j}$ and $H_{ij}$ are the first and second derivatives
 with respect to the back-ground metric $g_0$. We write $k(x,t):=\left(e^{-u}+\frac{2-\e}2 h'(r_\e^{1/2}(g)e^{(\e-2)u})
r_\e^{1/2}(g)e^{(\e-2)u}\right)$. We can check that $\forall (x,t)\in M\times[0,T_1]$ we have  $C_1>k(x,t)>C_2$.  Here the positive constants $C_1$ and $C_2$ are independent of $T$.

 From the positivity of $(F^{ij})$ and (\ref{54.1}), we have 
\begin{equation}
\label{5add1}\begin{array}{rcl}
0& \ge &  \vs\ds H_t- \frac{1}{2}\sum_{i,j}F^{ij} H_{ij} \\
&\ge &\vs\ds \frac 1 {2\s_2(g)} \tr \{T_{1}(S_g)\tilde \n_g^2 (e^{-u})\}+k(x,t) u_t
\\
&=&\vs\ds \frac12\sum_{i,j} F^{ij}
\{(e^{-u})_{ij}+u_i(e^{-u})_{j}+u_j(e^{-u})_{i}-u_l(e^{-u})_l\d_{ij}\}+k(x,t)u_t\\
&=&\vs\ds \frac12 e^{-u}\sum_{i,j} F^{ij}\{-u_{ij}-u_iu_j+|\n
u|^2\d_{ij}\}
+k(x,t) u_t\\
&=&\vs\ds \frac12 e^{-u}\sum_{i,j}  F^{ij}\{-w_{ij}+S(g_0)_{ij}+\frac
12|\n u|^2\d_{ij}\}
+k(x,t) u_t\\
&\ge & \vs\ds\frac12 e^{-u}\sum_{i,j} F^{ij}\{-w_{ij}+S(g_0)_{ij}\}+
k(x,t) u_t\\
&=&\ds\vs \frac12 e^{-u}\sum_{i,j} F^{ij}S(g_0)_{ij}+O(1)\log\s_2(g)+O(1) \\
&&- \frac 12k(x,t)(h(r_\e^{1/2}(g)e^{(\e-2)u})+s_\e(g)).
\end{array}\end{equation}
Here we have used $\sum_{i,j} F^{ij}w_{ij}=\sum_{i,j}\frac 1 {\s_2(g)}
\frac {\partial \s_2(g)}{\partial w_{ij} } w_{ij}=2$.
Since $g_0 \in \Gamma_2^+$, by Garding's inequality \cite{Garding},
\begin{equation}
\label{5add2}
\sum_{i,j}F^{ij}S(g_0)_{ij}=\sum_{i,j}\frac 1{\s_{2}(g)}\frac{\partial \s_{2}(g)}{w_{ij}} S(g_0)_{ij}
\ge 2e^{2u}\frac{\s^{1/2}_2(g_0)}{\s^{1/2}_2(g)}.\end{equation}
On the other hand, 
one can check $h(r_\e^{1/2}(g)e^{(\e-2)u})+s_\e(g)$ is bounded from above,
for $\|u\|_{C^2}$ is uniformly bounded.
Now from (\ref{5add1}) and (\ref{5add2}), we have
\[\begin{array}{rcl}
0 & \ge& \ds\vs e^{u}\frac{\s^{1/2}_2(g_0)}{\s^{1/2}_2(g)}+
O(1)\log\s_2(g)+O(1)\\
&\ge& \ds \frac{c_1}{\s^{1/2}_2(g)} +c_2\log \s_2(g)-c_3,\end{array}\]
for positive constants $c_1$, $c_2$ and $c_3$ independent of $T$.
Clearly, this inequality
implies that there is a constant $c_0>0$ independent of $T$
such that (\ref{54.1}) holds. Namely
\[\s_2(g) \ge c_0,\]
at point $(x_0, t_0)$. Hence, we have for any point $(x, t)\in M\times [0,T_1]$
\[\begin{array}{rcl}
&&\ds\log(e^{-2 u(x,t)}\s_2^{1/2}(g)(x,t))-\frac12 h(r_\e^{1/2}(g)e^{(\e-2)u(x,t)}) -e^{-u(x,t)}\\
& =& \ds\vs H(x,t)\\
&\ge& H(x_0,t_0)\\
&=& \ds\vs\log(e^{-2 u(x_0, t_0)}\s_2^{1/2}(g)(x_0, t_0)) -\frac12 h(r_\e^{1/2}(g)e^{(\e-2)u(x_0,t_0)})-e^{-u(x_0,t_0)}\\
&\ge& \log C_1-
e^{C},\end{array}\]
provided $e^{-2 u(x,t)}\s_2^{1/2}(g)(x,t)<1$. It follows that $\s_2(g)(x,t)
\ge C_2>0$ for some positive constant independent of $T$. \\
On $M\times [T_1, T]$, we  consider
a function $H:M\times [T_1,T]$ defined by 
$H=\log (e^{-\e u}\s_2^{1/2}(g))- e^{-u}$. By the  same argument,  
there is a constant $c>0$ independent of $T$ such that
$F(x,t)\ge c>0 \quad \hbox{ for any } (x,t)\in M\times [T_1,T]$. 
Hence, we deduce the claim (\ref{claim}).  This finishes the proof of the Proposition.
\qed

We remark that the Sobolev inequality, Theorem \ref{thmS0}, implies
that $r_\e$ has a positive lower bound. In the proof of Claim
(\ref{claim}), we avoided to use the Sobolev inequality.

\

\noindent{\it Proof of Proposition \ref{pro1}.} By  local estimates
established in \cite{GW1} (in fact a similar local estimates as in Theorem \ref{thmlocal1} hold), we can use the argument 
given in the proof of Proposition \ref{prop-e}  to show that the set of
solutions of (\ref{eq1-e}) with the bounded  ${\cal F}_2$ and $V_\e(e^{-2u}g_0)=1$ is compact for $\e\in (0,2]$. 
Hence, Proposition \ref{prop-e} implies that
$Y_\e$ is achieved by a function $u_\e$, which clearly is a solution of
(\ref{eq1-e}). We may assume that  $u_\e$ satisfies $V_\e (e^{-2u_\e}g_0)=1$
and
 \beq\label{eq-e2} \s_2(\uuu)=ce^{2(\e-2)u},\eeq 
 where $c=Y_\e$. For any fixed metric $g$, the function $\tilde {\cal F}_{2,\e}(g)$ is continuous on $\e$ so that $Y_\e$ is 
 semi-continuous 
 from above on $\e$. On the other hand, it follows from the H\"older's inequality, $Y_\e$ is 
 semi-continuous 
 from below on $\e$. Hence, $Y_\e$ is  continuous and we have
 \[\lim_{\e\to 0}Y_{\e}=Y_2(M,[g_0])<Y_2(\S^n).\]
If $\inf u_\e$ has a uniform lower bound, then the estimates established
in \cite{GW1} implies that $\|u_\e\|_{C^2}$ is uniformly bounded. By the
result of Evans-Krylov, $\|u_\e\|_{C^{2,\a}}$ is uniformly bounded for any $\a\in (0,1)$.
Hence $u_\e$, by taking a subsequence, converges strongly in $C^{2,\a}$
to $u_0$, which is a solution of (\ref{0.1}). Moreover, $u_0$ is a minimizer. 
Now suppose $\underline{\lim}_{\e \to 0} \inf u_\e=-\infty$. Let $(x_\e)\in M$ such that 
$u_\e(x_\e)=\min_{x\in M}u_\e(x)$. We consider a new function
\[ v_\e(y)=u(\exp_{x_\e} \d_\e y)-u_\e(x_\e)\]
and defined on $B_{\d_\e^{-1}}$ with a pull-back metric $g_\e:=(exp_{x_\e}\d_\e \cdot)^*g_0$,
where $\d_\e =e^{(1-\e/2)u_\e(x_\e)}$. Since $u_\e(x_\e)\to -\infty$,
$\d_\e\to 0$ as $\e\to 0$. And one can check that $B_{\d_\e^{-1}}$
tends to $\R^n$ and $g_\e$  to the standard Euclidean metric in any compact set in $\R^n$ for any $C^k$ norm. 
We can check that $v_\e$ satisfies 
the same Equation (\ref{eq-e2}) on $B_{\d_\e^{-1}}$ with  $S_{g_0}$ replaced by $S_{g_\e}$. By the local 
estimates in \cite{GW1}, $(v_\e)$ is uniformly bounded in $C^2$ on any fixed compact set. From 
the result of Evans-Krylov, it follows that
 $(v_\e)$ is uniformly bounded in $C^{2,\alpha}$  on any fixed compact set.
Hence,  $v_\e$ converges in any compact domain of $\R^n$
 to an entire solution $u$ of the following
equation on $\R^n$
\beq\label{l}
\s_2\left(\n^2 u+d u\otimes du-\frac 12 |\n u|^2 g_{\R^n}\right)= c_0e^{-4u},\eeq
with $c_0=Y_2(M, [g_0])$. It is easy to check that
 $\int_{\R^n} e^{-n u} dvol(g_{\R^n})\le 1$. In fact, we have
\[
\begin{array}{rcl}
\ds\int_{B_{\d_\e^{-1}}} e^{(2\e-n)v_\e}dvol({g_\e})&=&\ds \vs \d_\e^{-n}e^{(n-2\e)u_\e(x_\e)}\int_{B(x_\e,1)} e^{(2\e-n)u_\e}dvol({g_0})\\
&=&\ds \vs e^{(n/2-2)\e u_\e(x_\e)}\int_{B(x_\e,1)}
 e^{(2\e-n)u_\e}dvol({g_0})\leq V_\e (e^{-2u_\e}g_0)=1.
\end{array}
\]
Letting $\e\to 0$, the claim yields. By the classification of (\ref{l}) given in \cite{LiLi} or \cite{CGY3},
we have
$c_0 \ge Y_2(\S^n)>Y_2(M, [g_0])$, which contradicts $c_0=Y_2(M, [g_0])$. \qed

\section{Existence}

In this section, we will construct a conformal metric $\tilde g$ such that $\tilde {\cal F}_2(\tilde g)<Y_2(\S^n)$ and 
$\tilde g\in \Gamma_2^+$. Our construction 
is inspired from  Aubin's work \cite{Aubin}. See also \cite{Schoen}. The basic idea is 
to construct a suitable test functions. But the
 more delicate point in our case, as already mentioned in the Introduction,
is to keep 
the conformal metric in the admissible class $\Gamma_2^+$ as in \cite{GW3}. 
Fix a point $P\in M$. Assume $n\geq 5$. 
It follows from the work by Lee-Parker that there exists a conformal metric
$g_1$ on ${\cal M}$  such that in a normal coordinate system for $g_1$ at $P$
\beq
\label{curv1}
R=O(r^2),
\eeq
\beq
\label{curv2}
\D  R=-{\frac 16}|W(P)|^2,
\eeq
\beq
\label{curv3}
\mathrm{Ric}(P)=0,
\eeq
\beq
\label{curv4}
\sqrt{\det g_1}=1+O(r^5),
\eeq
where $r=|x|$. We denote
\beq
g_v=v^{-2}g_1,
\eeq
where
$$
v(x)=\left\{
\begin{array}{llll}
\lambda+ r^2,&\mbox{if }x\in B(0,r_0),\\
\lambda+ r^2_0,&\mbox{else. }
\end{array}
\right.
$$
We first need some estimates.
\begin{lem}
\label{lem4}
Assume
$$
A={g_1}^{-1}(\frac {\nabla^2_{g_1} v} v-{\frac 12}
{|\nabla_{g_1} v|^2\over v^2}{g_1}+S_{g_1}),
$$
where
$$
S_{g_1}={1\over n-2}(\mathrm{Ric}_{g_1}-{R\over 2(n-1)}{g_1}).
$$
Then we have
\beq
\label{s1}
\mathrm{tr}(A)={2n\lambda \over (\lambda +r^2)^2}+ { O(r^5)\over \lambda +r^2}+{R\over 2(n-1)}
\eeq
and
\beq
\label{s2}
\mathrm{tr}(A^2)= {4n \lambda^2\over (\lambda +r^2)^4}+
{2R\lambda\over (n-1)(\lambda +r^2)^2}-{\mathrm{Ric}(\nabla_{g_1} v, 
\nabla_{g_1} v) \over (\lambda +r^2)^2}+O(r).
\eeq
\end{lem}

\pr By definition, we have
\beq
\label{sigma}
\sigma_2(g_v)=v^4 \sigma_2({g_1}^{-1}({\nabla^2_{g_1} v\over v}-{1\over 2}{|\nabla_{g_1} v|^2\over v^2}{g_1}+S_{g_1})).
\eeq
It is clear that
$$
\begin{array}{llll}
\mathrm{tr}(A)&\vs\ds{={\D _{g_1} v\over v}-{n\over 2}{|\nabla_{g_1} v|^2\over v^2}
+\mathrm{tr}({g_1}^{-1}S_{g_1})}\\
&\ds{={\D _{g_1} v\over v}-{n\over 2}{|\nabla_{g_1} v|^2\over v^2}
+{R\over 2(n-1)}},
\end{array}
$$
where
$$
\D _{g_1} ={1\over\sqrt{\det {g_1}}}\sum_{i,j}{\partial\over\partial x^i}
(\sqrt{\det {g_1}}{g_1}^{ij}{\partial\over\partial x^j}).
$$
Recall
\beq
|\nabla_{g_1} v|^2=4r^2.
\eeq
In view of (\ref{curv4}), we have
\beq
\D _{g_1} v=2n+O(r^5).
\eeq
Therefore,
$$
\begin{array}{llll}
\vs \mathrm{tr}(A)&\ds{={ 2n+O(r^5)\over \lambda +r^2}-{n\over 2}{4r^2\over (\lambda +r^2)^2}
+{R\over 2(n-1)}}\\
&\ds{={2n\lambda \over (\lambda +r^2)^2}+ 
{ O(r^5)\over \lambda +r^2}+{R\over 2(n-1)}}.
\end{array}
$$
It is also easy to check that
$$
\begin{array}{llll}
\vs\mathrm{tr}(A^2)&=&\ds{{|\nabla^2_{g_1} v|^2\over v^2}+{n|\nabla_{g_1} v|^4\over 4v^4}+ \mathrm{tr}(({g_1}^{-1}S_{g_1})^2)}\\
&&\ds{-{|\nabla_{g_1} v|^2 \D _{g_1} v\over v^3}+
2\mathrm{tr}( {{g_1}^{-1}\nabla^2_{g_1} v \;{g_1}^{-1}S_{g_1}\over v})- 
{|\nabla_{g_1} v|^2 \over v} \mathrm{tr}({g_1}^{-1}S_{g_1})}.
\end{array}
$$
A simple calculation gives us 
\begin{eqnarray}
\label{second1} {n|\nabla_{g_1} v|^4\over
4v^4}&=& {4nr^4\over (\lambda +r^2)^4}, \\
\mathrm{tr}(({g_1}^{-1}S_{g_1})^2)&=&O(r^2), \\
-{|\nabla_{g_1} v|^2 \D _{g_1} v\over v^3}&=& -{8nr^2 \over (\lambda
+r^2)^3}+O(r), \\
 - {|\nabla_{g_1} v|^2 \over {v^2}}
\mathrm{tr}({g_1}^{-1}S_{g_1})&=& -\frac {2r^2R}{(n-1)(\lambda +r^2)^2}, \\
{g_1}^{-1}\nabla^2_{g_1} v &=& 2I + O(r^2) ,\\
\mathrm{tr}({g_1}^{-1}\nabla^2_{g_1} v
\;{g_1}^{-1}S_{g_1})&=& 2\mathrm{tr}({g_1}^{-1}S_{g_1})+
O(r^3)={R\over n-1}+ O(r^3) ,\\ \label{second2}
2\mathrm{tr}( {{g_1}^{-1}\nabla^2_{g_1} v \;{g_1}^{-1}S_{g_1}\over
v})&=& {2R\over (n-1)(\lambda +r^2)}+ O(r). \end{eqnarray}
To handle
$\ds{{|\nabla^2_{g_1} v|^2\over v^2}}$, we recall the Bochner's
formula \beq \langle\nabla(\D  v),\nabla v\rangle= -
|\nabla^2_{g_1} v|^2 +{1\over 2} \D (|\nabla_{g_1} v|^2)-
\mathrm{Ric}(\nabla_{g_1} v,\nabla_{g_1} v). \eeq 
Hence, we have
\beq \label{second3}
\begin{array}{llll}
\vs |\nabla^2_{g_1} v|^2&\ds{= - \langle\nabla(\D  v),\nabla v\rangle +
{1\over 2} \D (|\nabla_{g_1} v|^2)-
\mathrm{Ric}(\nabla_{g_1} v,\nabla_{g_1} v)}\\
\vs &\ds{=- \langle\nabla(2n+O(r^5)),\nabla v\rangle+ {1\over 2} \D (4r^2)-
\mathrm{Ric}(\nabla_{g_1} v,\nabla_{g_1} v)}\\
&=\ds{4n-
\mathrm{Ric}(\nabla_{g_1} v,\nabla_{g_1} v)+O(r^5).}
\end{array}
\eeq
>From (\ref{second1}), (\ref{second2}) and (\ref{second3}), we have
$$
\begin{array}{llll}
\vs \mathrm{tr}(A^2)&=&\ds{{4n\over (\lambda +r^2)^2}+ {4nr^4\over (\lambda +r^2)^4}- {8nr^2\over (\lambda +r^2)^3}+
{2R\over (n-1)(\lambda +r^2)}-{2Rr^2\over (n-1)(\lambda +r^2)^2}}\\
&&\ds{-{\mathrm{Ric}(\nabla_{g_1} v,\nabla_{g_1} v) \over (\lambda +r^2)^2}+O(r)}\\
&=&\ds{ {4n \lambda^2\over (\lambda +r^2)^4}+
{2R\lambda\over (n-1)(\lambda +r^2)^2}-
{\mathrm{Ric}(\nabla_{g_1} v,\nabla_{g_1} v) \over (\lambda +r^2)^2}+O(r).}
\end{array}
$$
\qed

\begin{lem}
\label{lem5}
Assume $\b\in ({1\over 2},{1\over 4})$. We have $\sigma_1(g_v)>0$ and $\sigma_2(g_v)>0$ in $B(0,\lambda^\b)$.
Moreover, if  $n\geq 9$, we have
\beq
\label{s3} \ba{rcl}
\ds\vs \int_{B(0,\lambda^\b)}\sigma_2(g_v) dvol({g_v})&=&
\lambda^{-{n\over 2}+2}\{2n(n-1)B+C\D  R(0)\lambda^2\\
&&\ds +
O\left(\lambda^{5\over 2}+\lambda^{n({1\over 2}-\b)}+
\lambda^{2+(n-8)({1\over 2}-\b)}\right)\}\ea
\eeq
and
\beq
\label{s4}
\int_{B(0,\lambda^\b)}dvol({g_v})= 
\lambda^{-{n\over 2}}
\left[B+O\left(\lambda^{5/2}+\lambda^{n({1/2}-\b)}\right)\right],
\eeq
where the constants $B,C$ are given by
\beq
B=\int_{\R^n}{1\over (1+|x|^2)^n}dx,
\eeq
\beq
C=\int_{\R^n}\left({|x|^2\over 2n(1 +|x|^2)^{n-2}}+{2|x|^4\over
n(n+2)(1 +|x|^2)^{n-2}}\right)dx>0.
\eeq
\end{lem}

\pr It follows directly from (\ref{s1}) and (\ref{s2})
$$
\sigma_1(g_v)= v^2\left({ 2n\lambda\over (\lambda +r^2)^2}+\frac{R}{2(n-1)}
+O(r^3)\right),
$$
\beq\label{add1} \sigma_2(g_v)={v^4\over 2}\left[{4n
(n-1)\lambda^2\over (\lambda +r^2)^4}+ \frac {2\l
R}{(\l+r^2)^2}+{\mathrm{Ric}(\nabla_{g_1} v,\nabla_{g_1} v) \over
(\lambda +r^2)^2}+O(r)\right]. \eeq Thus, the first part of lemma
is clear. On the other hand, we obtain
\beq\label{add2}
\begin{array}{llll}
&\vs \ds{\int_{B(0,\lambda^\b)}\sigma_2(g_v) dvol({g_v})}\\
=&\vs\ds \int_{B(0,\lambda^\b)}{1\over (\lambda +r^2)^n}\{2n(n-1)\lambda^2 +
R\lambda(\lambda +r^2)^2 \\
&\ds + \frac 12\mathrm{Ric}((\nabla_{g_1} v,\nabla_{g_1} v)
(\lambda +r^2)^2+O(r(\lambda +r^2)^4)\} (1+O(r^5))dx.
\end{array}
\eeq
We can calculate
$$
\begin{array}{lll}
\vs \mathrm{Ric}(\nabla_{g_1} v,\nabla_{g_1} v)
 &=&\ds 4\sum_{i,j}R_{ij}(x)x^ix^j \\
\vs &=&\ds 4\sum_{i,j}( R_{ij}(0)+  \sum_{k}R_{ij,k}(0)x^k+\sum_{k,l}{1\over 2}R_{ij,kl}(0)x^kx^l)x^ix^j +O(r^5)\\
&=&\ds 2\sum_{i,j,k,l}R_{ij,kl}(0)x^kx^l x^ix^j +O(r^5).
\end{array}
$$
It is known that (see \cite{Aubin})
$$
{1\over r^4w_{n-1}}\int_{S(r)} \sum_{i,j,k,l}R_{ij,kl}(0)x^kx^l x^ix^j d\O=
{2\over n(n+2)}\D  R(0)
$$
and
$$
R(x)={1\over 2}\sum_{i,j}R_{,ij}(0)x^ix^j +O(r^3),
$$
 where $S(r)$ is the geodesic sphere of radius equal to $r$, $d\O$ is the volume element on the unit sphere $\S^{n-1}\subset\R^n$ 
and $w_{n-1}$ is the volume of the unit sphere $\S^{n-1}$.
Therefore,
$$
\begin{array}{llll}
\vs &\ds{\int_{B(0,\lambda^\b)}\sigma_2(g_v) dvol({g_v})}\\
\vs =&\ds{\int_{B(0,\lambda^\b)} {2n(n-1)\lambda^2\over (\lambda +r^2)^n}+{\lambda r^2 \D  R(0)\over
2n(\lambda +r^2)^{n-2}}dx  }\\
\vs &\ds{+\int_{B(0,\lambda^\b)}\left( {2r^4 \D  R(0)\over
n(n+2)(\lambda +r^2)^{n-2}}
+{O(r)\over (\lambda +r^2)^{n-4}}\right)dx}\\
=&\ds{\lambda^{-n/2+2}\int_{B(0,{\lambda}^{\b-1/2})} \left({2n(n-1)\over (1 +r^2)^n} +a(x)\lambda^2\D
R(0)\right)dx+O(\lambda^{-n/2+4+1/2})},
\end{array}
$$
where
$$
a(x)={|x|^2\over 2n(1 +|x|^2)^{n-2}}+{2|x|^4\over n(n+2)(1
+|x|^2)^{n-2}}.
$$
Thus, (\ref{s3}) yields. Similarly, we can estimate
$$
\begin{array}{lll}
\vs \ds{\int_{B(0,\lambda^\b)}dvol({g_v})}
\vs&=&\ds{\int_{B(0,\lambda^\b)} v^{-n} \sqrt{\det g_1} dx}\\
\vs&=&\ds{\int_{B(0,\lambda^\b)}{1+O(r^5) \over (\lambda +r^2)^n}dx}\\
\vs&=&\ds{\lambda^{-n/2}\int_{B(0,{\lambda}^{\b-1/2})} {dx\over (1 +r^2)^n} +O(\lambda^{-n/2+5/2})}\\
&=&\ds{\lambda^{-n/2}\left[B+O\left(\lambda^{5/2}+
\lambda^{n({1/2}-\b)}\right)\right]}.
\end{array}
$$
Therefore, we finish the proof.
\qed

\begin{lem}\label{lem6} Let $g_1$ as above and $\gamma\in (0,2)$ be given. Assume $n\geq 9$. For sufficiently small
$\d>0$ such that $\lambda^{1/4}>>\d>>\lambda^{1/2}$, there exists a constant $1>\d_1>\d$ and a function $u:
B_{\d_1}\to \R$ satisfying :
\begin{itemize}
\item[(0)]$\d_1^{n-4\over 2}=(\frac 2\gamma-1)\l^{-1}\d^{n\over 2}(1+o(1))$,
\item[(1)] The metric $\tilde g=e^{-2u}g_1$ has positive
$\Gamma_2$-curvatures,
\item[(2)] $u=\log (\l+|x|^2)+b_0$ for $|x|\le \d$,
\item[(3)] $u=\gamma\log|x|$ for $|x|\ge \d_1$,
\item[(4)] $ vol(B_{\d_1}\backslash B_{\d}, \tilde g)\le C\ds\left(\frac{
\d^{\frac{n+4-n\gamma}{2(2-\gamma)}}}{\lambda}\right)^{2n(2-\gamma)/(n-4)}$,
\item[(5)] $\int_{B_{\d_1}\backslash B_{\d}}\s_2(\tilde g)
d vol(\tilde g)\le C \d^{4+n(1-\gamma)}\l^{-3+2\gamma}  $,
\end{itemize} where $b_0$ satisfies (\ref{late}) below.
\end{lem}
\pr We want to find a function $u$ with $u'(r)=\frac {\a(r)} r$. The
Schouten tensor of $\tilde g=e^{-2u}g_1$ is
\begin{equation}
\begin{array}{rcl} \ds\vs S(\tilde g)_{ij}&=& \ds \n^2_{ij} u+\n_i u\n_j
u-\frac {|\n
u|^2} 2 {(g_1)}_{ij}+S(g_1)_{ij}\\
&=&\ds \frac {2\a}{2r^2} \d_{ij}-\frac {\a^2}{2r^2} (g_1)_{ij} +\left (\frac {\a'} r+\frac
{\a^2-2\a}{r^2}\right) \frac{x_ix_j}{r^2}+S(g_1)_{ij}+O(r^2)\frac
{\a} {r^2},
\end{array}
\end{equation}
so that
\begin{equation}
\begin{array}{rcl} \ds\vs S(\tilde g)^i_{ j}&=&\ds \frac {2\a-\a^2}{2r^2} \d_{ij} +\left (\frac {\a'} r+\frac
{\a^2-2\a}{r^2}\right) \frac{x_ix_j}{r^2}+S(g_1)^i_{ j}+O(r^2)\frac
{\a} {r^2},
\end{array}
\end{equation}
since it follows from Gauss Lemma that $\sum_i (g_1)^{ij}x_i=x_j$. We look for a function $\a(r)\in (\gamma,2)$ for all $r\in (\d,\d_1)$.
Hence one can find  a fixed constant $A>0$ independent of $\l$ such that
\begin{equation}
  S(\tilde g)^i_{ j}\ge \frac {2\a-\a^2-Ar^2\a}{2r^2} \d_{ij} +\left (\frac {\a'} r+\frac
{\a^2-2\a}{r^2}\right) \frac{x_ix_j}{r^2}
\end{equation}
and
\begin{equation}
S(\tilde g)^i_{ j}\le \frac {2\a-\a^2+Ar^2\a}{2r^2} \d_{ij} +\left (\frac {\a'} r+\frac
{\a^2-2\a}{r^2}\right) \frac{x_ix_j}{r^2}.
\end{equation}
Consequently, we obtain
\[\s_2(\tilde g)
>e^{4 u}\frac{(n-1)}{2}\left(\frac {2\a-\a^2-Ar^2\a}{2r^2}\right)^2
\left(n-4+4\frac {r\a'-Ar^2\a}{2\a-\a^2-Ar^2\a}\right)\]
and
\[\s_1(\tilde g)
>e^{2 u}\left(\frac {2\a-\a^2-Ar^2\a}{2r^2}\right)
\left(n-2+2\frac {r\a'-Ar^2\a}{2\a-\a^2-Ar^2\a}\right).\]
We want to
find an $\a$ satisfying
\[\a=\left\{
\begin{array}{ll}
\ds\vs \frac {2r^2}{\l+r^2}, & \hbox{ if } |x|\le \d,\\
\ds\vs \mbox{solution of }(\ref{e1}), & \hbox{ if } |x|\in (\d,\d_1),\\ \gamma,&\hbox{ if }
|x|\ge \d_1.
\end{array}\right.\]
Such a function can be found as follows.  First we solve the
following equation
\begin{equation}\label{e1}
\frac {n-4} 4+\frac{
r\a'-Ar^2\a}{2\a-\a^2-Ar^2\a}=0.\end{equation}
Recall  this is
the Bernoulli differential equation. One can find a
general solution of (\ref{e1}) as follows.
\[\frac 1{\a}=r^{n-4\over 2} e^{-\frac {n A r^2} 2}\left(\int_1^r \frac{4-n}
4 \frac 1{t^{n-2\over 2}} e^{\frac {nA t^2} 2 } dt+c\right).\] Set
\[H(r)=-\frac {An} 2\int_1^r \frac 1{t^{n-6\over 2}} e^{\frac {nAt^2} 2
} dt.\]
We have
\[
\begin{array}{rcl}
\a&=&\ds\vs \frac 2{1+2 a_1 r^{n-4\over 2} e^{-\frac {nA
r^2}2 }+2H(r) r^{n-4\over 2} e^{-\frac {nA
r^2}2 }}\\
&=&\ds\frac2{1+2a_1 r^{n-4\over 2}+2G(r)},\end{array}\]
where \[G= a_1
r^{n-4\over 2}(e^{-\frac {nAr^2} 2 }-1)+H(r)  r^{n-4\over 2} e^{-\frac {nA
r^2}2 }.\] Here the constant $a_1$ is determined by
\[ \a(\d)= \frac{2\d^2}{\l+\d^2}.\]
We have the estimate
\begin{equation}\label{e2}a_1=\frac {\l}{2\d^{n\over 2}}(1+o(1)),
\end{equation}
since we use the fact $\l^{1/4}>>\d$. Define $\d_1$ by $\a(\d_1)=\gamma$. We have
\begin{equation}\label{e4}
\d_1^{n-4\over 2}=\left(\frac2\gamma-1\right)\l^{-1}\d^{n\over 2}(1+o(1))
\end{equation}
so that $1>>\d_1>>\d$. Note that $n\geq 9$. Hence, for all $r\in (\d,\d_1)$ we have
\begin{equation}\label{e3}
G(r)=O(1) r^2,
\end{equation}
so that for all $r\in(\d,\d_1)$
\begin{equation}\label{e4.1}u(r)=\frac 4{4-n}\log
(r^{4-n\over 2}+2a_1)+a_2,\end{equation} where \[a_2=(\gamma-2)\log
\d_1-\frac 4{4-n}\log \frac 2\gamma +o(1).\]
For $r<\d$ we have \[u(r)=\log(\l+r^2)+b_0,\] where
\begin{equation}\label{late}
b_0=(\gamma-2)\log\d_1+O(1),\end{equation}
where we use $\d>>\l^{1/2}$.
In view of (\ref{e1}), we have
\begin{equation}\label{f1}
  (2\a-\a^2+Ar^2\a)\left(n-4+4\frac
{r\a'+Ar^2\a}{2\a-\a^2+Ar^2\a}\right)=2nAr^2\a=O(1)r^2.
\end{equation}
We also have for all $r\in (\d,\d_1)$
\begin{equation}
  \a(r)\in (\gamma,2)
\end{equation}
and
\begin{equation}\label{f2}
  2-\a+Ar^2= \frac {4a_1 r^{n-4\over 2}}{1+2a_1r^{n-4\over 2}}+O(r^2).
\end{equation}

 Now we can check that
\[\begin{array}{l}
\ds\vs\quad\quad  \int_{B_{\d_1}\backslash B_\d}\s_2(\tilde g)
dvol(\tilde g)\\ \le \ds C(n)\int_{B_{\d_1}\backslash B_\d}
e^{(4-n)u}\left(\frac {2\a-\a^2+Ar^2\a}{2r^2}\right)^2
\left(n-4+4\frac {r\a'+Ar^2\a}{2\a-\a^2+Ar^2\a}\right) dvol(g_1)\\
\le \ds\vs O(1)\int_{\d}^{\d_1} e^{(4-n)a_2
}(r^{4-n\over 2}+2a_1)^4\frac {1}{r^2} (2-\a+Ar^2) r^{n-1} dr\\
\le \ds\vs O(1)\int_{\d}^{\d_1} \d_1^{(n-4)(2-\gamma)}
r^{5-n} {a_1 r^{n-4\over 2}}{(1+2a_1 r^{n-4\over 2})^4}dr\\
\le \ds\vs O(1)\int_{\d}^{\d_1} \d_1^{(n-4)(2-\gamma)}  {a_1 r^{3-{n\over 2}} dr }\\
\le \ds\vs O(1) \d_1^{(n-4)(2-\gamma)}\d^{4-{n\over 2}} a_1=O(1)
\d^{n(1-\gamma)+4}\l^{-3+2\gamma}
\end{array}\]
and
\[\begin{array}{ll}
\ds\vs\quad\quad vol(B_{\d_1}\backslash B_{\d}, \tilde g)& = \ds \int_{B_{\d_1}\backslash B_\d}
e^{-nu} dvol(g_1)\\
&\le \ds\vs O(1)\int_{\d}^{\d_1} e^{-na_2
}(r^{4-n\over 2}+2a_1)^{4n/(n-4)} r^{n-1} dr\\
&= \ds\vs O(1)\int_{\d}^{\d_1} \d_1^{n(2-\gamma)} r^{-1-n} (1+2a_1 r^{n-4\over 2})^{4n/(n-4)}dr\\
&\le \ds\vs O(1) \d_1^{n(2-\gamma)}\d^{-n} =O(1)\left(\frac{
\d^{\frac{n+4-n\gamma}{2(2-\gamma)}}}{\lambda}\right)^{2n(2-\gamma)/(n-4)}.
\end{array}\]
Therefore, after smoothing $u$, we get a desired $u$.\qed

We write $g_0=e^{-2u_0}g_1$. In the following result, we try to connect the initial metric $g_0$ to some tube object.
More precisely, we prove the following lemma.\\

\begin{lem}\label{lem7} Let $g_0\in \Gamma_2^+$ and the geodesic ball $B(0,r_0)$ as above. Assume that $n\geq 5$.
For any given $\gamma\in (0,2)$, then there is a
conformal metric $\tilde g= e^{-2u}g_1$ of positive $\Gamma_2$-curvatures on $
B(0,r_0)\setminus \{0\}$ 
satisfying :
\begin{itemize}
\item[(1)] The metric $\tilde g=e^{-2u}g_1$ has positive
positive $\Gamma_2$-curvatures;
\item[(2)] $u=\gamma\log|x|$ for $|x|\le r_2$;
\item[(3)] $u=u_0(x)+b_1$ for $|x|\ge r_1$;
\end{itemize} where $r_2<r_1<r_0$ and $b_1$ is a constant.
\end{lem}
\pr
We write $u(x)=w(r) + \xi(r) u_0(x)$ where $\xi(r)$ is some  cut-off function 
equals to 1 near of $r_0$ and to 0 near 0,
and $w$ with $w'(r)=\frac {\a(r)} r$, where $\a$ is equel to 0 near $r_0$.
As before, the
Schouten tensor of $\tilde g=e^{-2u}g_1$ is
\begin{equation}
\begin{array}{rcl} \ds\vs S(\tilde g)_{ij}&=& \ds \n^2_{ij} w+\n_i w\n_j
w + \n_i w\n_j
 (\xi u_0)+ \n_i  (\xi u_0)\n_j
w\\
&&\ds -\left(\frac {|\n
w|^2} 2 +\langle \n w,\n (\xi u_0)\rangle \right)(g_1)_{ij}+S(e^{-2\xi u_0}g_1)_{ij}\\
\end{array}
\end{equation}
so that
\begin{equation}
\begin{array}{rcl} \ds\vs S(\tilde g)^i_{ j}&=&\ds \frac {2\a-\a^2}{2r^2} \d_{ij} +\left (\frac {\a'} r+\frac
{\a^2-2\a}{r^2}\right) \frac{x_ix_j}{r^2}+S(e^{-2\xi u_0}g_1)^i_{ j}+O(r+|\n (\xi u_0)|)\frac
{\a} {r}.
\end{array}
\end{equation}
Fix $\e\in (0,{2-\gamma\over 5})$ and let $C_1$ bound the term  $O(r+|\n   u_0|)$.
Set $r_4=\min({r_0\over2},{1\over 2},{\e\over 2(1+C_1)})$. For some small $r_3$ to be fixed later,
we want to $\a$ decrease from $\gamma$ to 0 in $(r_3, r_4)$ and $\xi\equiv 1$ in $(r_3, r_0)$.
In $B_{r_0}\setminus B_{r_3}$, we write $A=S(\tilde g)-S( g_0)$. Therefore
$$
\s_2(\tilde g)=e^{4(w+u_0)}\s_2(A+ S( g_0)).
$$
We want $A+ S( g_0)\in \Gamma_2^+$ in $B_{r_0}\setminus B_{r_3}$. It is clear in $B_{r_4}\setminus B_{r_3}$
$$
A\ge \left(\frac {2\a-\a^2-\e\a}{2r^2} \d_{ij} +\left (\frac {\a'} r+\frac
{\a^2-2\a}{r^2}\right) \frac{x_ix_j}{r^2}\right).
$$
This gives
\[\s_2(A)
>e^{4 u}\frac{(n-1)}{2}\left(\frac {2\a-\a^2-\e\a}{2r^2}\right)^2
\left(n-4+4\frac {r\a'-\e\a}{2\a-\a^2-\e\a}\right)\]
and
\[\s_1(A)
>e^{2 u}\left(\frac {2\a-\a^2-\e\a}{2r^2}\right)
\left(n-2+2\frac {r\a'-\e\a}{2\a-\a^2-\e\a}\right).
\]
We see that for all small $\d>0$,
\begin{equation}
\label{in1}
\a(r)=\frac{(2-5\e)\d}{\d+r^{{1\over 2}-{5\over 4}\e}}
\end{equation}
solves the equation
\begin{equation}
{1\over 4}(2\a-\a^2-\e\a)=-r\a'+\e\a.
\end{equation}
We choose some $r_5<r_4$ and a non increasing function $\a$ in $(r_5,r_4)$ such that $\a(r_4)=0$, $\a(r_5)>0$
and $\tilde g\in \Gamma_2^+$ in  $B_{r_4}\setminus B_{r_5}$ by openness of $ \Gamma_2^+$. Now we choose a suitable $\d$
in (\ref{in1}) and take a small $r_6<r_5$ such that $\a(r_6)=\gamma$. From the above calculations, we see that $A\in \Gamma_2^+$ in $B_{r_5}\setminus B_{r_6}$ so that $ S( \tilde g)\in \Gamma_2^+$ in $B_{r_5}\setminus B_{r_6}$. Now we set $\a(r)=\gamma$ for all $r<r_6$ and $r_3=r_6$. We see that
there exists some cut-off function $\xi$ such  that $\xi(r)=1\forall r>r_7$,   $\xi(r)=0\forall r<r_8$ and $r^2S(e^{-2\xi u_0}g_1)$
is small in $B_{r_7}\setminus B_{r_8}$ where $r_8<r_7<r_6$. Thus we can choose such suitable cut-off function such that $\tilde g$
in $\Gamma_2^+$. Now it is sufficient  to choose some $r_2<r_8$ and $r_1=r_4$. Finally, we obtain the desired $u$ by smoothing it.
\qed

The construction of such  metrics of positive $\Gamma_2$-curvatures is motivated by the method introduced
by Gromov-Lawson \cite{GromovLawson} in their study of metrics of 
positive scalar curvature.
See also for the constructions of other positive metrics in \cite{MiWang} and
\cite{GLW}.
 Now we can prove the main result in this section.\\

\begin{thm}
\label{thmY} Let $(M,g_0)$ be a compact, oriented
 Riemannian manifold with $\s_2(g_0)>0$. Assume that $n\geq 9$. Then there exists $\tilde g\in [g_0]$ such that
 \begin{equation}
  \label{eq1-thm4.1}
  \tilde g\in \Gamma_2^+
 \end{equation}
 and
 \begin{equation}
 \label{eq2-thm4.1}
 \tilde {\cal F}_2(\tilde g)<Y_2(\S^n).
 \end{equation}
 \end{thm}
\pr We fix somme $\gamma\in (1,2)$ and let the geodesic ball $B(0,r_0)$ w.r.t. $g_1$
 as above. We define a conformal metric $\tilde g$ as follows. Let $r_2<r_1<r_0$ as
 in Lemma \ref{lem7} and 
 set $\d=\lambda^{\beta}$ with $\beta\in ({1\over 4},{1\over 2})$ for any small
 $\lambda$. Find $\d_1$ as in Lemma \ref{lem6}. Now for any small 
 $\l$ with $\d_1<r_2$, define $\tilde g$  on $B_{\d_1}$ by Lemma \ref{lem6}
 and on $B_{r_0}\backslash B_{r_2}$ by Lemma \ref{lem7}. 
 And on $M\backslash B(0,r_0)$, $\tilde g=e^{-2b_1}g_0$, where the constant
$b_1$ is given in Lemma \ref{lem7}. Since on $B_{r_2}\backslash B_{\d_1}$
the metrics constructed in Lemma \ref{lem6} and Lemma \ref{lem7} are the same,
$\tilde g$ is smooth.
 From Lemmas \ref{lem6} and \ref{lem7}, we know
 (\ref{eq1-thm4.1}) holds. In the following, we keep the notations
 of the geodesic ball with respect to the background metric $g_1$. 
 By the Lemmas \ref{lem5}, \ref{lem6} and \ref{lem7}, we can estimate
\begin{equation}
\label{eq3-thm4.1}
\int_{B_{\d_1}\backslash B_{\d}}\s_2(\tilde g)
d vol(\tilde g)\le C \d_1^{(n-4)(2-\gamma)}\d^{4-{n\over 2}} a_1,
\end{equation}

\begin{equation}
\label{eq4-thm4.1}
\int_{M\backslash B_{\d_1}}\s_2(\tilde g)
d vol(\tilde g)\le C \d_1^{(n-4)(2-\gamma)}\d_1^{4-n},
\end{equation}

\begin{equation}
\label{eq5-thm4.1}
\begin{array}{llll}
\ds \int_{ B_{\d}}\s_2(\tilde g)
d vol(\tilde g)&=& \d_1^{(n-4)(2-\gamma)} \lambda^{-{n\over 2}+2}\left[2n(n-1)B+C\D  R(0)\lambda^2
\right.\\
&&
+\left.O\left(\lambda^{5\over 2}+\lambda^{n({1\over 2}-\b)} 
+\lambda^{2+(n-8)({1\over 2}-\b)}\right)\right],
\end{array}
\end{equation}

\begin{equation}
\label{eq6-thm4.1}
vol(M, \tilde g)\ge vol( B_{\d}, \tilde g)= \d_1^{n(2-\gamma)}
\lambda^{-{n\over 2}}\left[B+O\left(\lambda^{5/2}+\lambda^{n({1/2}-\b)}\right)
\right].
\end{equation}
We choose some $\beta\in ({1\over 4}, {n-4\over 2n})$ so that we obtain
$$
\d^{4-{n\over 2}} a_1=o(\lambda^{-{n\over 2}+4})\,\,\mbox{and}\,\,\d_1^{4-n}
=o(\lambda^{-{n\over 2}+4}).
$$
As a consequence, we get
$$
\tilde {\cal F}_2(\tilde g)\leq B^{4-n\over n} 
[2n(n-1)B+C\D  R(0)\lambda^2 +o(\lambda^2)].
$$
Recall $2n(n-1) B^{4\over n}=Y_2(\S^n)$ and $\D  R(0)<0$. Therefore, we deduce (\ref{eq2-thm4.1}) provided $\lambda$ is sufficiently small. Hence,
we finish the proof.
\qed

\

\noindent{\it Acknowledgement.} A part of the work was carried out while
the first author was visiting Max Planck Institute for Mathematics in the
 Sciences,
Leipzig and while the second author
was visiting Universit\'e Paris XII-Val de Marne.
We  would like to thank both institutions and J\"urgen Jost and Frank Pacard 
for warm hospitality.
The second author also would like to thank Pengfei Guan for his many helpful
discussions. The main result of the paper was announced at a talk given by 
the second author in Magdeburg University, Dec. 3, 2004. When this work was done, 
we learned that Weimin Sheng, Neil Trudinger and Xujia Wang 
had announced a similar result.  Their result  can now 
be found in \cite{STW}, in which they showed  that (3) is true for $n>4$.
Hence, the $\s_2$-Yamabe problem is completely solved.


\

%

\end{document}